\documentclass{article}
\usepackage{amsthm, amssymb, latexsym, amsmath}
\usepackage[bookmarks,bookmarksnumbered]{hyperref}
\newtheorem{theorem}{Theorem}[section]
\newtheorem{lemma}[theorem]{Lemma}
\newtheorem{proposition}[theorem]{Proposition}
\newtheorem{corollary}[theorem]{Corollary}

\newtheorem{remark}{Remark}

\newtheorem{definition}{Definition}
\newenvironment{them}{\textbf{Theorem}\it}{}
\include{now_macro}
\huge
\begin{document}
\title{Perturbation Expansion and $N$-th Order Fermi Golden Rule of the Nonlinear Schr\"odinger Equations}
\author{Zhou Gang\thanks{This paper is part of the author's
Ph.D thesis.} \footnote{Supported by NSERC under Grant NA7901.} }
\maketitle \centerline{\small{Department of Mathematics, University
of Toronto, Toronto, Canada}}
%\fixNumberingInArticle
\setlength{\leftmargin}{.1in}
\setlength{\rightmargin}{.1in}%\NowFootNum %\fixNumberingInArticle
\normalsize \vskip.1in \setcounter{page}{1}
\setlength{\leftmargin}{.1in} \setlength{\rightmargin}{.1in}
\section*{Abstract}
In this paper we consider generalized nonlinear Schr\"odinger
equations with external potentials. we compute the forth and the
sixed order Fermi Golden Rules (FGR), conjectured in ~\cite{GS2,
GS3}, which is used in a study of the asymptotic dynamics of trapped
solitons.
\section{Introduction}\label{Sec:intro}
It is well known that the eigenvalues embedded in a continuous
spectrum of a self-adjoint operator, $L$, are unstable under generic
perturbations. Intuitively, one can think of time periodic solutions
(`bound states') of the equation, $i\dot u=L u$, corresponding to
such eigenvalues, leaking their energy to the continuous spectrum
solutions (`scattering states') which in turn is radiated to
infinity. The coupling between time-periodic and continuous spectrum
solutions was computed in the second order of perturbation theory by
P. A. M. Dirac in 1927 (see~\cite{Di}, resulting in an elegant
expression. If this expression is nonzero, then the decay mentioned
above takes place. E. Fermi called this condition ``Golden Rule No.
2". In physics literature it is known as the Fermi Golden Rule
(FGR). The FGR, computed first for atomic and molecular states,
appears in many other areas of physics e.g. in non-equilibrium
statistical mechanics. Of course, due to energy conservation the
isolated eigenvalues are stable under reasonable time independent
perturbations.

The Fermi Golden Rule was introduced to nonlinear Hamiltonian PDE in
~\cite{S} where it was used to prove that time-periodic solutions
for linear (and nonlinear) wave equations are unstable under generic
nonlinear perturbations. It was used extensively in investigating long-time dynamics of
solitary waves in the nonlinear Schr\"odinger
~\cite{BP2,BuSu,SW1,SW2,SW3,TY1,TY2,TY3} and wave \cite{SW4} equations.

However, in nonlinear problems with a small parameter, $h$, of the
quasi-classical nature, the second order expression obtained by
Dirac vanishes and, in principle, one can hope to get non-vanishing
expression only in the $2N$th order where $N=O(\frac{1}{h})$. (If
one thinks about the FGR as a condition that an 'atom is ionized by
a single photon', in this case because the ionization energy is much
larger than the photon frequency, this 'ionization' is due to a
multi-photon process.) This type of situation occurs in particular
in the problem of the relaxation of a solution to the ground state
in generalized nonlinear Schr\"odinger equations with external
slowly varying potentials (Gross-Pitaevskii equations) described
below.

In this paper we find, for the first time, the 'nonlinear' forth and
sixth order FGRs for the relaxation problem mentioned above. As was
conjectured in ~\cite{GS2,GS3}, they are of the same form as the second
order FGR and as, presumably, the higher order FGRs. To obtain the
forth and sixth order FGRs we have to perform additional normal-form-type
transformations in the perturbation theory we developed in
~\cite{GS2,GS3}.

\textbf{Problem.} Now we describe the problem in more details. We
consider the generalized nonlinear Schr\"odinger equation (NLS) in
dimension $d\neq 2$ with an external potential $V_{h}:\
\mathbb{R}^{d}\rightarrow \mathbb{R}$,
\begin{equation}\label{NLS}
i\frac{\partial\psi}{\partial t}=-\Delta
\psi+V_{h}\psi-f(|\psi|^{2})\psi.
\end{equation}
Here $h>0$ is a small parameter giving the length scale of the
external potential in relation to the length scale of the $V_{h}=0$
solitons (see below), $\Delta$ is the Laplace operator and $f(s)$ is
a nonlinearity to be specified later. We normalize $f(0)=0$. Such
equations arise in the theory of Bose-Einstein condensation
\footnote[1]{In this case Equation (~\ref{NLS}) is called the
Gross-Pitaevskii equation.}, nonlinear optics, theory of water waves
\footnote[2]{In the last two areas $V_{h}$ arises if one takes into
account impurities and/or variations in geometry of the medium and
is, in general, time-dependent.}and in other areas.

To fix ideas we assume the potentials to be of the form
$V_{h}(x):=V(hx)$ with $V$ smooth and decaying at $\infty.$ Thus for
$h=0,$ Equation (~\ref{NLS}) becomes the standard generalized
nonlinear Schr\"odinger equation (gNLS)
\begin{equation}\label{gNLS}
i\frac{\partial \psi}{\partial t}=-\Delta
\psi+\mu\psi-f(|\psi|^{2})\psi,
\end{equation} where $\mu:=V(0).$
For a certain class of nonlinearities, $f(|\psi|^{2})$, (see Section
~\ref{ExSt}), there is an interval $\mathcal{I}_{0}\subset
\mathbb{R}$ such that for any $\lambda\in \mathcal{I}_{0}$ Equation
(~\ref{gNLS}) has solutions of the form
$e^{i(\lambda-\mu)t}\phi_{0}^{\lambda}(x)$ where
$\phi_{0}^{\lambda}\in \mathcal{H}_{2}(\mathbb{R}^{n})$ and
$\phi_{0}^{\lambda}>0.$ Such solutions (in general without the
restriction $\phi_{0}^{\lambda}>0$) are called the \textit{solitary
waves} or \textit{solitons} or, to emphasize the property
$\phi_{0}^{\lambda}>0,$ the \textit{ground states}. For brevity we
will use the term \textit{soliton} applying it also to the function
$\phi_{0}^{\lambda}$ without the phase factor $e^{i(\lambda-\mu)t}.$

Equation (~\ref{gNLS}) is translationally and gauge invariant. Hence
if $e^{i(\lambda-\mu)t}\phi_{0}^{\lambda}(x)$ is a solution for
Equation (~\ref{gNLS}), then so is $$
e^{i(\lambda-\mu)t}e^{i\alpha}\phi_{0}^{\lambda}(x+a),\ \text{for
any}\ a\in \mathbb{R}^{n},\text{and}\ \alpha\in [0,2\pi).$$ This
situation changes dramatically when the potential $V_{h}$ is turned
on. In general, as was shown in ~\cite{Floer,Oh1,ABC} out of the
$(n+2)$-parameter family
$e^{i(\lambda-\mu)t}e^{i\alpha}\phi_{0}^{\lambda}(x+a)$ only a
discrete set of two-parameter families of solutions to Equation
(~\ref{NLS}) bifurcate: $e^{i\lambda
t}e^{i\alpha}\phi^{\lambda}(x),$ $\alpha\in [0,2\pi)$ and
$\lambda\in \mathcal{I}$ for some $\mathcal{I}\subseteq
\mathcal{I}_{0}$, with $\phi^{\lambda}\equiv \phi_{h}^{\lambda}\in
\mathcal{H}_{2}(\mathbb{R}^{n})$ and $\phi^{\lambda}>0$. Each such
family centers near a different critical point of the potential
$V_{h}(x).$ It was shown in ~\cite{Oh2} that the solutions
corresponding to minima of $V_{h}(x)$ are orbitally (Lyapunov)
stable and to maxima, orbitally unstable. We call the solitary wave
solutions described above which correspond to the minima of
$V_{h}(x)$ \textit{trapped solitons} or just \textit{solitons} of
Equation (~\ref{NLS}) omitting the last qualifier if it is clear
which equation we are dealing with.

A basic question about soliton solutions is whether they are
asymptotically stable, i.e. whether for initial condition of
(~\ref{NLS}) sufficiently close to a trapped soliton  $\{e^{i\gamma}
\phi^{\lambda}(x)\}$ the solution converges (relaxes) in a local
norm, up to a phase factor, to another trapped soliton,
$$\psi(x,t)-e^{i\gamma(t)}\phi^{\lambda_{\infty}}(x)\rightarrow 0.$$

We observe
that (~\ref{NLS}) is a Hamiltonian system with conserved energy (see
Section ~\ref{HaGWP}) and, though orbital (Lyapunov) stability is
expected, the asymptotic stability is a subtle matter. To have
asymptotic stability the system should be able to dispose of excess
of its energy, in our case, by radiating it to infinity. The
infinite dimensionality of a Hamiltonian system in question plays a
crucial role here. This phenomenon as well as a general class of
classical and quantum relaxation problems was pointed out by J.
Fr\"ohlich and T. Spencer [Private Communication].

Another important property is their effective dynamics. Namely, one would
like to show that if an initial condition is close (in the weighted
norm
$\|u\|_{\nu,1}:=\|(1+|x|^{2})^{\frac{\nu}{2}}u\|_{\mathcal{H}^{1}}$
for sufficiently large $\nu$) to the soliton
$e^{i\gamma_{0}}\phi^{\lambda_{0}},$ with $\gamma_{0}\in \mathbb{R}$
and $\lambda_{0}\in \mathcal{I}$ ($\mathcal{I}$ as above), then the
solution, $\psi(t),$ of Equation (~\ref{NLS}) can be written as
\begin{equation}\label{decom1}
\psi(x,t)=e^{i\gamma(t)}[e^{ip(t)\cdot
x}\phi^{\lambda(t)}(x-a(t))+R(x,t)],
\end{equation}
where
%$\phi_{hz_1z_2}^{\lambda}(x):= \phi^{\lambda}(x-a)$ with $a:=(z_1, z_2)$,
%$\gamma(t)\in \mathbb{R},$
$\|R(t)\|_{-\nu,1}\rightarrow 0$, $\lambda(t)\rightarrow
\lambda_{\infty}$ for some $\lambda_{\infty}$ as $t\rightarrow
\infty$ and the soliton center $a(t)$ and momentum $p(t)$ evolve
according to certain effective equations of motion.

\textbf{Results}.As in ~\cite{GS2,GS3} we assume that either $d=1$
and the potential is even or $d>2$ and the potential is spherically
symmetric and the initial condition symmetric with respect to
permutations of the coordinates. In this case the soliton relaxes to
the ground state along the radial direction. This limits the number
of technical difficulties we have to deal with. We expect that our
techniques extend to the general case when the soliton spirals
toward its equilibrium.

It is shown in ~\cite{GS2, GS3}, under certain conditions, that the
ground state is asymptotically stable and the effective equations
for the parameters $a$ and $p$ are close to Newton's equations. It
was conjectured in ~\cite{GS2,GS3} that one of the conditions (which
plays a key role) is equivalent to an explicit $2N$th order FGR.
This conjectured is true for $N=1$ due to the works
~\cite{BP2,BuSu,SW1,SW2,SW3,TY1,TY2,TY3}. In this paper we prove
this conjecture for $N=2$ and $3.$ In other words we establish the
FGR for nonlinear problem in the $4$th and $6$th orders.

Moreover, we find more precise effective equations for $a$ and $p.$
We show that there exist a function $z(t)$ and a constant
$Z_{N+1,N}\in \mathbb{C}$ satisfying
$z(t)=p(t)+ia(t)+O(|p|^{2}+|a|^{2})$ and
$$\frac{1}{2}\frac{d}{dt}|z|^{2}=ReZ_{N+1,N}|z|^{2N+2}+O(|z|^{2N+3})$$ with
$N=O(1/h)$ being an integer depends on $\lambda$ and $h.$ For
$N=2,3$ we find the explicit form of $ReZ_{N+1,N}$ which in
particular shows that it is always non-positive and negative
generically. The explicit form is given by $2N$th order FGR
mentioned above.

In the present paper by using normal forms we simplify the
expressions for $z$ and the other parameters, thus make the proof
more transparent than ~\cite{GS2,GS3}.

\textbf{Previous results}. We refer to ~\cite{GS1} for a detailed
review of the related literature. Here we only mention results of
~\cite{Cu,Buslaev,BP2,BuSu,SW1,SW2,SW3,TY1,TY2,TY3} which deal with
a similar problem.  Like our work, ~\cite{SW1,SW2,SW3,TY1,TY2,TY3}
study the ground state of the NLS with a potential. However, these
papers deal with the near-linear regime in which the nonlinear
ground state is a bifurcation of the ground state for the
corresponding Schr\"{o}dinger operator $-\Delta+V(x).$ The present
paper covers highly nonlinear regime in which the ground state is
produced by the nonlinearity (our analysis simplifies considerably
in the near-linear case).

Papers ~\cite{Cu,Buslaev,BP2,BuSu} consider the NLS without a
potential so the corresponding solitons, which were described above,
are affected only by a perturbation of the initial conditions which
disperses with time leaving them free. In our case they, in
addition, are under the influence of the potential and they relax to
an equilibrium state near a local minimum of the potential.

As was mentioned above, the relaxation of the solution to
the ground state (the asymptotic stability of the ground state) was
shown in ~\cite{GS2,GS3} under some assumptions on the potential
$V$, nonlinearity $f$ and the linearized operator $L(\lambda)$,
which in particular include these related to the FGR which was
explained above.

The paper is organized in the following way. In Section
~\ref{Sec:intro} we introduce the concept of Fermi Golden Rules and
outline its applications and importance in proving the asymptotic
stability of trapped solitons of nonlinear Schr\"odinger equations.
In Section ~\ref{HaGWP} we show the conservation laws and the local
well-posedness of (~\ref{NLS}). In Section ~\ref{ExSt} we formulate
the conditions on the nonlinearity $f$ and the potential $V$ so that
there exists a soliton manifold to (~\ref{NLS}) and it is stable. In
Section ~\ref{subslinear} we linearize the solution to (~\ref{NLS})
around the soliton get a linear operator. Moreover we analyze its
spectrum and put some assumptions. In Section ~\ref{MainTHM} we
state the main theorem. In Section ~\ref{SEC:effective} we separate
the 'useless' parts from the equations of parameters, to prepare for
the computation of forth ($N=2$) and sixth ($N=3$) order Fermi
Golden Rules which will be proved in Section ~\ref{SEC:ProofFGR} and
~\ref{sec:CaseN3} respectively.

\textbf{Notation}. As customary we often denote derivatives by
subindices as in
$\phi^{\lambda}_{\lambda}=\frac{\partial}{\partial\lambda}\phi^{\lambda}$
for $\phi^{\lambda}=\phi^{\lambda}(x).$ However, the subindex $h$
signifies always the dependence on the parameter $h$ and not the
derivatives in $h.$ The Sobolev and $L^{2}$ spaces are denoted by
$\mathcal{H}^{k}$ and $\mathcal{L}^{2}$ respectively.
%\section*{Acknowledgment}

\textbf{Acknowledgment}. This paper is part of the author's Ph.D
thesis requirement. The author wishes to thank his advisor,
I.M.Sigal, for introducing the problem and useful discussions.

We are grateful to  J. Colliander, S. Cuccagna, S. Dejak, J.
Fr\"ohlich, Z. Hu, W. Schlag, A. Soffer, G. Zhang, V. Vougalter and,
especially, V.S. Buslaev for fruitful discussions. Upon the
completion of this paper we learn from T.P.Tsai that he also
obtained the explicit form of Fermi Golden rule for $N=2$ [Private
Communications].
\section{Hamiltonian Structure and GWP}\label{HaGWP}
Equation (~\ref{NLS}) is a Hamiltonian system on Sobolev space
$\mathcal{H}^{1}(\mathbb{R},\mathbb{C})$ viewed as a real space
$\mathcal{H}^{1}(\mathbb{R},\mathbb{R})\oplus
\mathcal{H}^{1}(\mathbb{R},\mathbb{R})$ with the inner product
$(\psi,\phi)=Re\int_{\mathbb{R}}\bar{\psi}\phi$ and with the
simpletic form
$\omega(\psi,\phi)=Im\int_{\mathbb{R}}\bar{\psi}\phi.$ The
Hamiltonian functional is: $$H(\psi):=\int
[\frac{1}{2}(|\psi_{x}|^{2}+V_{h}|\psi|^{2})-F(|\psi|^{2})],$$ where
$F(u):=\frac{1}{2}\int_{0}^{u}f(\xi)d\xi.$

Equation (~\ref{NLS}) has the time-translational and gauge
symmetries which imply the following conservation laws: for any
$t\geq 0,$ we have
\begin{enumerate}
 \item[(CE)] conservation of energy: $H(\psi(t))=H(\psi(0));$
 \item[(CP)]
 conservation of the number of particles: $N(\psi(t))=N(\psi(0)),$ where $N(\psi):=\int
 |\psi|^{2}.$
\end{enumerate}
To address the global well-posedness of (~\ref{NLS}) we need the
following condition on the nonlinearity $f$. Below, $s_{+}=s$ if
$s>0$ and $=0$ if $s\leq  0$.
\begin{enumerate}
 \item[(fA)] The nonlinearity $f$ satisfies the
 estimate $|f^{'}(\xi)|\leq c(1+|\xi|^{\alpha-1})$ for some $\alpha\in
 [0,\frac{2}{(d-2)_{+}})$  and $|f(\xi)|\leq c(1+|\xi|^{\beta})$ for some $\beta\in[0,\frac{2}{d}).$
\end{enumerate}

The following result can be found in ~\cite{Cazenave}.\\
\begin{them}
Assume that the nonlinearity $f$ satisfies the condition (fA), and
that the potential $V$ is bounded. Then Equation (~\ref{NLS}) is
globally well posed in $\mathcal{H}^{1}$, i.e. the Cauchy problem
for Equation (~\ref{NLS}) with a datum $\psi(0)\in \mathcal{H}^{1}$
has a unique solution $\psi(t)$ in the space $\mathcal{H}^{1}$ and
this solution depends continuously on $\psi(0)$. Moreover $\psi(t)$
satisfies the conservation laws (CE) and (CP).
\end{them}
\section{Existence and Orbital Stability of Solitons}\label{ExSt}
In this section we review the question of existence of the solitons
(ground states) for Equation (~\ref{NLS}). Assume the nonlinearity
$f:\mathbb{R}\rightarrow \mathbb{R}$ is smooth and satisfies
\begin{enumerate}
 \item[(fB)] There exists an interval $\mathcal{I}_{0}\in \mathbb{R}^{+}$
s.t. for any $\lambda\in\mathcal{I}_{0}$, $-\infty\leq
\displaystyle\overline{\lim}_{s\rightarrow +
 \infty}\frac{f(s)}{s^{\frac{2}{d-2}}}\leq 0$
and $\frac{1}{\xi}\int_{0}^{\xi}f(s)d
 s>\lambda$ for some constant $\xi$, for $d>2$;
and
 $$U(\phi,\lambda):=-\lambda\phi^{2}+\int_{0}^{\phi^{2}}f(\xi)d\xi$$
 has a smallest positive root $\phi_{0}(\lambda)$ such that
 $U_{\phi}(\phi_{0}(\lambda),\lambda)\not=0$, for $d=1$.
 \end{enumerate}

It is shown in ~\cite{BL, Str} that under Condition (fB) there
exists a spherical symmetric positive solution $\phi^{\lambda}$ to
the equation
\begin{equation}\label{soliton}
-\Delta\phi^{\lambda}+\lambda\phi^{\lambda}-f((\phi^{\lambda})^{2})\phi^{\lambda}=0.
\end{equation}
\begin{remark}
Existence of soliton functions $\phi^{\lambda}$ for $d=2$ is proved
in ~\cite{Str} under different conditions on $f$.
\end{remark}
When the potential $V$ is present, then some of the solitons above
bifurcate into solitons for Equation (~\ref{NLS}). Namely, let, in
addition, $f$ satisfy the condition $|f^{'}(\xi)|\leq
c(1+|\xi|^{p}),$ for some $p< \infty$, and $V$ satisfy the condition
\begin{enumerate}
\item[(VA)] $V$ is smooth and $0$ is a non-degenerate local
minimum of $V$.
\end{enumerate}
Then, similarly as in ~\cite{Floer,Oh1} one can show that if $h$ is
sufficiently small, then for any $\lambda\in \mathcal{I}_{0V}$,
where
$$\mathcal{I}_{0V}:=\{\lambda|\lambda>-\displaystyle\inf_{x\in\mathbb{R}}\{V(x)\}\}\cap\{\lambda|\lambda+V(0)\in
\mathcal{I}_{0}\},$$ there exists a unique soliton
$\phi^{\lambda}\equiv\phi_{h}^{\lambda}$ (i. e. $\phi^{\lambda}\in
\mathcal{H}_{2}(\mathbb{R})$ and $\phi^{\lambda}>0$) satisfying the
equation
$$-\Delta \phi^{\lambda}+(\lambda+V_{h})\phi^{\lambda}-
f((\phi^{\lambda})^{2})\phi^{\lambda}=0$$ and the estimate
$\phi^{\lambda}=\phi^{\lambda+V(0)}_{0}+O(h^{3/2})$ where
$\phi_{0}^{\lambda}$ is the soliton of Equation (~\ref{soliton}).

Let $\delta(\lambda):=\frac{1}{2}\|\phi^{\lambda}\|_{2}^{2}$. It is
shown in ~\cite{GSS1} that the soliton $\phi^{\lambda}$ is a
minimizer of the energy functional $H(\psi)$ for a fixed number of
particles $N(\psi)=constant$ if and only if
%\begin{equation}\label{Stab}
$\delta^{'}(\lambda)>0.$
%\end{equation}
%(\textbf{GZ:  Does one need additional conditions of} $f$?)
Moreover, it is shown in ~\cite{We2,GSS1} that under the latter
condition the solitary wave $\phi^{\lambda}e^{i\lambda t}$ is
orbitally stable.
%(\textbf{GZ:  Does one need additional conditions
%of} $f$?)
Under more restrictive conditions (see ~\cite{GSS1}) on $f$ one can
show that the open set
%there exists an interval
%$\mathcal{I}$ such that if $\lambda\in \mathcal{I}$ then
%$\delta^{'}(\lambda)>0.$
%
%In what follows we set
\begin{equation}
\mathcal{I}:=\{\lambda\in \mathcal{I}_{0V}:\delta'(\lambda)>0\}
\end{equation}
is non-empty. Instead of formulating these conditions we assume in
what follows that the open set $\mathcal{I}$ is non-empty and
$\lambda\in \mathcal{I}$.

Using the equation for $\phi^{\lambda}$ one can show that if the
potential $V$ is redially symmetric then there exist constants $c,\
\delta>0$ such that
%\begin{equation}\label{expondecay}
$|\phi^{\lambda}(x)|\leq ce^{-\delta|x|}\ \text{and}\
|\frac{d}{d\lambda}\phi^{\lambda}|\leq ce^{-\delta|x|},$
%\end{equation}
and similarly for the derivatives of $\phi^{\lambda}$ and
$\frac{d}{d\lambda}\phi^{\lambda}$.

%%%%%%%%%%%%%%%%%%%%%%%%%%%%%%%%%%%%%%%%%%%%%%%%%%%%%%%%%%%%%%%%%%%%%%%
%%%%%%%%%%%%%%%%%%%%%%%%%%%%%%%%%%%%%%%%%%%%%%%%%%%%%%%%%%%%%%%%%%%%%%%%%
\section{Linearized Equation and Resonances}\label{subslinear} We
rewrite Equation (~\ref{NLS}) as $\frac{d\psi}{dt}=G(\psi)$ where
the nonlinear map $G(\psi)$ is defined by
%\begin{equation}\label{symmetryG}
$G(\psi)=-i(-\Delta+\lambda+V_{h})\psi+if(|\psi|^{2})\psi.$
%\end{equation}
%appearing on the right hand side of Equation (~\ref{NLS}).
Then the linearization of Equation (~\ref{NLS}) can be written as
$\frac{\partial\chi}{\partial t}=\partial G(\phi^{\lambda})\chi$
where $\partial G(\phi^{\lambda})$ is the Fr\'echet derivative of
$G(\psi)$ at $\phi$. It is computed to be
\begin{equation}\label{defineoperator}
\partial
G(\phi^{\lambda})\chi=-i(-\Delta+\lambda+V_{h})\chi+if((\phi^{\lambda})^{2})\chi+2if^{'}((\phi^{\lambda})^{2})(\phi^{\lambda})^{2}Re\chi.
\end{equation}
This is a real linear but not complex linear operator. To convert it
to a linear operator we pass from complex functions to real
vector-functions $\chi\longleftrightarrow \vec{\chi}=\left(
\begin{array}{lll}
\chi_{1}\\
\chi_{2}
\end{array}
\right), $ where $\chi_{1}=Re\chi$ and $\chi_{2}=Im\chi.$ Then
$\partial G(\phi^{\lambda})\chi\longleftrightarrow
L(\lambda)\vec{\chi}$ where the operator $L(\lambda)$ is given by
\begin{equation}\label{operaL}
L(\lambda) :=  \left(
\begin{array}{lll}
0&L_{-}(\lambda)\\
-L_{+}(\lambda)&0
\end{array}
\right),
\end{equation}
with
%\begin{equation}\label{firstoperator}
$L_{-}(\lambda):=-\Delta+V_{h}+\lambda-f((\phi^{\lambda})^{2}),$
%\end{equation}
and
%\begin{equation}\label{secondoperator}
$L_{+}(\lambda):=-\Delta+V_{h}+\lambda-f((\phi^{\lambda})^{2})-2f^{'}((\phi^{\lambda})^{2})
(\phi^{\lambda})^{2}.$
%\end{equation}
The operator $L(\lambda)$ is extended to the complex space
$\mathcal{H}^{2}(\mathbb{R},\mathbb{C})\oplus
\mathcal{H}^{2}(\mathbb{R},\mathbb{C}).$ If the potential $V_{h}$ in
Equation (~\ref{NLS}) decays at $\infty$, then by a general result
%(see e.g. ~\cite{HS,RSIV}),
$$\sigma_{ess}(L(\lambda))=(-i\infty,-i\lambda]\cap
[i\lambda,i\infty).$$

The eigenfunctions of $L(\lambda)$ are described in the following
theorem (cf ~\cite{GS1}, ~\cite{GS2}).

%%%%%%%%%%%%%%%%%%%%%%%%%%%%%%%%%%%%%%%%%%%%%%%%%%%%%%%%%%%%%%%%%%%%%%%%%%%%%%%%%%%%%%%%%
\begin{theorem}\label{mainpo}
Let $V$ satisfy Condition (VA) and $|h|$ be sufficiently small. Then
the operator $L(\lambda)$ has at least $2d +2$ eigenvectors and
associated eigenvectors with eigenvalues near zero: two-dimensional
space with the eigenvalue 0 and a $2d $-dimensional space with
non-zero imaginary eigenvalues $\pm i\epsilon_j(\lambda),$
$$\epsilon_j(\lambda):=h\sqrt{2e_j} +o(h),$$ where $e_j$ are
eigenvalues of the Hessian matrix of $V$ at value $x=0,$
$V^{''}(0)$. The corresponding eigenfunctions
$\left(\begin{array}{lll}
\xi_j\\
\pm i\eta_j
\end{array}
\right)$ are related by complex conjugation and satisfy
$$
\xi_j= \sqrt{2}\displaystyle \partial_{x_{k}}\phi_{0}^{\lambda}
 +o(h)\
\text{and}\ \eta_j= -h \sqrt{e_j}\displaystyle
x_{j}\phi^{\lambda}_{0} +o(h),$$ and $\xi_{i}$ and $\eta_{j}$ are
real.
\end{theorem}

%%%%%%%%%%%%%%%%%%%%%%%%%%%%%%%%%%%%%%%%%%%%%%%%%%%%%%%%%%%%%%%%%%%%%%%
\begin{remark}
The zero eigenvector $\left(
\begin{array}{lll}
0\\
\phi^{\lambda}
\end{array}
\right)$ and the associated zero eigenvector $\left(
\begin{array}{lll}
\partial_{\lambda}\phi^{\lambda}\\
0
\end{array}
\right)$ are related to the gauge symmetry $\psi(x,t)\rightarrow
e^{i\alpha}\psi(x,t)$ of the original equation and the $2d$
eigenvectors $\left(\begin{array}{lll}
\xi_j\\
\pm i\eta_j
\end{array}
\right)$ with $O(h)$ eigenvalues originate from the zero
eigenvectors $\left(
\begin{array}{lll}
\partial_{x_{k}}\phi_{0}^{\lambda}\\
0
\end{array}
\right), k=1,2,\cdot\cdot\cdot,d,$ and the associated zero
eigenvectors $\left(
\begin{array}{lll}
0\\
x_{k}\phi_{0}^{\lambda}
\end{array}
\right),\ k=1,2,\cdot\cdot\cdot,d,$ of the $V=0$ equation due to the
translational symmetry and to the boost transformation
$\psi(x,t)\rightarrow e^{ib\cdot x}\psi(x,t)$ (coming from the
Galilean symmetry), respectively.
\end{remark}

For  $d\geq 2$ we will be interested in permutationally symmetric
functions, $g\in \mathcal{L}^{2}(\mathbb{R}^{d})$, characterized as
$$g(x)=g(\sigma x)\ \text{for any}\ \sigma\in S_{d}$$ with $S_{d}$ being the group of
permutation of $d$ indices and $\sigma
(x_{1},x_{2},\cdot\cdot\cdot,x_{d}):=(x_{\sigma(1)},x_{\sigma(2)},
\cdot\cdot\cdot,x_{\sigma(d)}).$
\begin{remark}\label{remark2}
For any function of the form $e^{ip\cdot x}\phi(|x-a|)$ with $a
\parallel p$, there exists a rotation $\tau$ such that
the function $e^{ip\cdot \tau
{x}}\phi(|\tau{x}-a|)=e^{i\tau^{-1}p\cdot x}\phi(|x-\tau^{-1}a|)$ is
permutationally symmetric. Such families describe wave packets with
the momenta directed toward or away from the origin.
\end{remark} If for $d\geq 2$ the potential $V(x)$ is spherically
symmetric, then $V^{''}(0)=\frac{1}{d}\Delta V(0)\cdot Id_{n\times
n}$, and therefore the eigenvalues $e_{j}$ of $V^{''}(0)$ are all
equal to $ \frac{1}{d}\Delta V(0)$. Thus we have
%%%%%%%%%%%%%%%%%%%%%%%%%%%%%%%%%%%%%%%%%%%%%%%%%%%%%%%%%%%%%%%%%%%%%%%%%%%%%
%%%%%%%%%%%%%%%%%%%%%%%%%%%%%%%%%%%%%%%%%%%%%%%%%%%%%%%%%%%%%%%%%%%%%%%%%%%%%%
\begin{corollary}\label{mainpo1}
Let $d\geq 2$ and $V$ satisfy Condition (VA) and let $V$ be
spherically symmetric. Then $L(\lambda)$ restricted to permutational
symmetric functions has $4$ eigenvectors or associated eigenvectors
near zero: two-dimensional space with eigenvalue 0; and
two-dimensional space with the non-zero imaginary eigenvalues $\pm
i\epsilon(\lambda)$, where $$\epsilon(\lambda)=h\sqrt{\frac{2\Delta
V(0)}{d}}+o(h),$$ and with the eigenfunctions
$\left(\begin{array}{lll}
\xi(\lambda)\\
\pm i\eta(\lambda)
\end{array}
\right)$, where $\xi$ and $\eta$ are real, and permutation symmetric
functions satisfying
$$\xi(\lambda)=\sqrt{2}\displaystyle\sum_{n=1}^{d}\frac{d}{dx_{n}}\phi_{0}^{\lambda}
+O(h)\ \text{and}\ \eta(\lambda)=-h\sqrt{\frac{1}{d}\Delta
V(0)}\displaystyle\sum_{n=1}^{d}x_{n}\phi^{\lambda}_{0} +O(h^{3/2})
.$$
\end{corollary}

Besides eigenvalues, the operator $L(\lambda)$ may have resonances
at the tips, $\pm i\lambda$, of its essential spectrum (those tips
are called thresholds). Recall the notation $\alpha_+ := \alpha$ if $\alpha>0$ and $=0$ of
$\alpha\leq 0$.

%To define the resonance we write the operator $L(\lambda)$ as
%$L(\lambda)=L_{0}(\lambda)+V_{big}(\lambda),$

%%%%%%%%%%%%%%%%%%%%%%%%%%%%%%%%%%%%%%%%%%%%%%%%%%%%%%%%%%%%%%%%%%%%%
%%%%%%%%%%%%%%%%%%%%%%%%%%%%%%%%%%%%%%%%%%%%%%%%%%%%%%%%%%%%%%%%%%%%
\begin{definition}
Let $d\neq 2$. A function $h$ is called a resonance function of
$L(\lambda)$ at $\mu=\pm i\lambda$ if $h\not\in \mathcal{L}^{2}$,
$|h(x)|\leq c\langle x\rangle^{-(d-2)_+}$ and $h$ is $C^{2}$ and
solves the equation
$$(L(\lambda)-\mu)h=0.$$
\end{definition}

Note that this definition implies that for $d>2$ the resonance
function $h$ solves the equation $(1+K(\lambda))h=0$ where
$K(\lambda)$ is a family of compact operators given by $K(\lambda)
:= (L_{0}(\lambda)-\mu+0)^{-1}V_{big}(\lambda)$. Here
%\begin{equation}\label{Lo}
$L_{0}(\lambda):=\left(
\begin{array}{lll}
0&-\Delta+\lambda\\
\Delta-\lambda&0
\end{array}
\right)$
%\end{equation}
and
\begin{equation}\label{Vbig}
V_{big}(\lambda):=\left(
\begin{array}{lll}
0&V_{h}-f((\phi^{\lambda})^{2})\\
-V_{h}+f((\phi^{\lambda})^{2})+2f^{'}((\phi^{\lambda})^{2})(\phi^{\lambda})^{2}&0
\end{array}
\right).
\end{equation}

In this paper we make the following assumptions on the point
spectrum and resonances of the operator $L(\lambda):$
\begin{enumerate}
 \item[(SA)] $L(\lambda)$ has only $4$ standard and associated
eigenvectors in the permutation symmetric subspace.
 \item[(SB)] $L(\lambda)$ has no resonances at $\pm i\lambda$.
\end{enumerate}

The discussion and results concerning these conditions, given in
~\cite{GS1}, suggested strongly that Condition (SA) is satisfied for
a large class of nonlinearities and potentials and Condition (SB) is
satisfied generically. Elsewhere
%In ~\cite{GSV}
we show this using earlier results of ~\cite{CP, CPV}. We also
assume the following condition
\begin{enumerate}
 \item[(FGR)] Let $N$ be the smallest positive integer
such that $\epsilon(\lambda)(N+1)> \lambda,\ \forall \lambda \in I$.
Then $Re Z_{N+1,N}< 0$ where $Z_{m,n},\ m,n=1,2,\cdots,$ are the
functions of $V$ and $\lambda,$ defined in Equation (~\ref{FGRZ})
below (see also (~\ref{FGR2})).
\end{enumerate}

We expect that Condition (FGR) holds generically. Theorem
~\ref{GOLD:maintheorem} below shows that $Re Z_{n+1,n}= 0$ if $n<N.$

We expect the following is true:\\
(a) if for some $N_{1}(\geq N),$ $ReZ_{n+1,n}=0$ for $n<N_{1},$ then
$ReZ_{N_{1}+1,N_{1}}\leq 0$ and (b) for generic
potentials/nonlinearities there exists an $N_{1}(\geq N)$ such that
$ReZ_{N_{1}+1,N_{1}}\not=0.$ Thus Condition (FGR) could have been
generalized by assuming that $Re Z_{N_1 +1,N_{1}}< 0$ for some $N_1
\geq N$ such that $ReZ_{n+1,n}=0$ for $n<N_{1}$. We took $N=N_1$ in
order not to complicate the exposition.

The following form of $Re Z_{N+1,N}$
 \begin{equation}\label{FGRY}
 Re Z_{N+1,N}=Im\langle \sigma_{1}(L(\lambda)-(N+1)i\epsilon(\lambda)-0)^{-1}F,F\rangle\leq 0
 \end{equation}
 for some function
 $F$ depending on $\lambda$ and $V$ and $\sigma_{1}:=\left(
 \begin{array}{lll}
 0&-1\\
 1&0
 \end{array}
 \right)$, is proved in ~\cite{BuSu, TY1, TY2, TY3, SW4} for $N=1,$ and in the present paper for $N=2,3$.
 We conjecture that this formula holds for any $N$. [See more explicit forms of $F$ for $N=1,2,3$ in Theorem ~\ref{GOLD:maintheorem}
 below].

 Condition (FGR) is related to the Fermi Golden Rule condition which
 appears whenever time-(quasi)periodic, spatially localized solutions
 become coupled to radiation. In the standard case it says that this coupling is effective
 in the second order ($N=1$) of the perturbation theory and therefore it leads to instability of such
 solutions. In our case these time-periodic solutions are
 stationary solutions
$$c_1 \left(
\begin{array}{lll}
\xi\\
 i\eta
\end{array}
\right)e^{i\epsilon(\lambda) t} +c_2 \left(
\begin{array}{lll}
\xi\\
- i\eta
\end{array}
\right)e^{-i\epsilon(\lambda) t}$$
 of the linearized equation $\frac{\partial\vec{\chi}}{\partial
 t}=L(\lambda)\vec{\chi}$ and the coupling is realized through the
 nonlinearity. Since the radiation in our case is "massive"$-$ the
 essential spectrum of $L(\lambda)$ has the gap
 $(-i\lambda,i\lambda)$, $\lambda>0,$ $-$ the coupling occurs only in
 the $N-$th order of perturbation theory where $N$ is given in
 Condition (FGR).

The rigorous form of the Fermi Golden Rule for the linear
Schr\"odinger  equations was introduced in ~\cite{BS}. For nonlinear
waves and Schr\"odinger equations the Fermi Golden Rule and the
corresponding condition were introduced in ~\cite{S} and, in the
present context, in ~\cite{SW4, BuSu, BP2, TY1, TY2, TY3}.

Recall that a function $g\in \mathcal{L}^{2}(\mathbb{R}^{d})$ is
permutational symmetric if
$$g(x)=g(\sigma x)\ \text{for any}\ \sigma\in S_{d}$$ with $S_{n}$ being the group of
permutation of $d$ indices and
$$\sigma
(x_{1},x_{2},\cdot\cdot\cdot,x_{d}):=(x_{\sigma(1)},x_{\sigma(2)},\cdot\cdot\cdot,x_{\sigma(d)}).$$
\section{Main Theorem}\label{MainTHM}
The following are two main theorems of this paper.
\begin{theorem}\label{GOLD:maintheorem}
Let $h$ be sufficiently small, the potential $V$ be spherically
symmetric and the initial date $\psi(0)$ be permutational symmetric
if $d\geq 3$. If $f$ satisfies the conditions (fA)-(fC) and (GWP);
$V$ satisfies the conditions (VA)-(VC) and (GWP); $\lambda\in
\mathcal{I}$, and the spectral conditions (SA) and (SB) are
satisfied, then there exists $c,\epsilon_{0}>0$ such that, if
\begin{equation}\label{InitCond}\inf_{\gamma\in
\mathbb{R}}\{\|\psi_0-e^{i\gamma}(\phi^{\lambda}+z_{1}^{(0)}\xi+iz_{2}^{0}\eta)\|_{\mathcal{H}^{k}}+
\|(1+x^{2})^{\nu}[\psi_0-e^{i\gamma}(\phi^{\lambda}+z_{1}^{(0)}\xi+iz_{2}^{0}\eta)]\|_{2}\}\leq
c|(z_{1}^{0},z_{2}^{0})|^{2}
\end{equation} with
$|(z_{1}^{0},z_{2}^{0})|\leq \epsilon_{0}$ and $z_{n}^{0} \ n=1,2$
being real, some large constant $\nu>0$ and with $k=[\frac{d}{2}]+3$
if $d\geq 3,$ and $k=1$ if $d=1,$ then there exist smooth functions
$\gamma,\ z_{1},\ z_{2}:\mathbb{R}^{+}\rightarrow \mathbb{R},$
$\lambda:\mathbb{R}^{+}\rightarrow \mathcal{I}$ and
$R:\mathbb{R}^{3}\times\mathbb{R}^{+}\rightarrow \mathbb{C}$ such
that
\begin{equation}\label{Decom}
\psi(x,t)=e^{i\int_{0}^{t}\lambda(s)ds}e^{i\gamma(t)}[\phi^{\lambda}+p_{1}\phi_{\lambda}^{\lambda}+ip_{2}\phi^{\lambda}+(z_{1}+p_{3})\xi+i(z_{2}+p_{4})\eta+R]\end{equation}
where $p_{n}(z,\bar{z}): \mathbb{C}\rightarrow \mathbb{R},\
n=1,2,3,4,$ are polynomials of $z:=z_{1}+iz_{2}$ and $\bar{z}$ and
of order $|z|^{2}$, the function $R(t)$ can be decomposed as
$$
\vec{R}=\displaystyle\sum_{2\leq m+n\leq k\leq
2N}R_{mn}(\lambda)z^{m}\bar{z}^{n}+R_{k}
$$ where for some large constant $\nu$ $\langle x\rangle^{-v}R_{m,n}\in \mathcal{L}^{2}$, $\|\langle x\rangle^{-v}R_{2N}\|_{2}\lesssim (1+t)^{-\frac{2N+1}{2N}}.$

The functions $\lambda$, $\gamma,$ $z$ have the following properties
\begin{enumerate}
\item[(A)] There exists a $\lambda_{\infty}\in \mathcal{I}$ such that $|\lambda(t)-\lambda_{\infty}|\leq
c(1+t)^{-\frac{1}{2N}}$ as $t\rightarrow\infty$ and moreover
\begin{equation}\label{ExpanLambda}
\dot\lambda=\sum_{2\leq m+n\leq
2N+1}\Lambda_{m,n}z^{m}\bar{z}^{n}+Remainder
\end{equation} where
 $\Lambda_{m,n}=0$ if $m,n\leq N+1;$
\item[(B)] \begin{equation}\label{ExpanGamma}
\dot\gamma=\sum_{2\leq m+n\leq
2N+1}\Gamma_{m,n}z^{m}\bar{z}^{n}+Remainder\end{equation} where
$\Gamma_{m,n}=0$ if $m,n\leq N+1$ and $m\not=n;$ $\Gamma_{m,m}$ is
real if $m\leq N$;
\item[(C)] $|z|\leq c(1+t)^{-\frac{1}{2N}}$ and moreover
\begin{equation}\label{FGRZ}
\dot{z}=-i\epsilon(\lambda)z+\sum_{2\leq m+n\leq
2N+1}Z_{m,n}z^{m}\bar{z}^{n}+Remainder
\end{equation} $Z_{n+1,n}(\lambda)$ is purely imaginary if $n< N;$ $Z_{m,n}(\lambda)=0$ if $m,n\leq N+1$ and $m\not=n+1;$ moreover when
$N=1,2,3$ we have
\begin{equation}\label{FGR}
Re Z_{N+1,N}=Im\langle
\sigma_{1}(L(\lambda)-(N+1)i\epsilon(\lambda)-0)^{-1}F,F\rangle\leq
0
\end{equation} where $F$ is some vector function
such that $(L(\lambda)-i(N+1)\epsilon(\lambda)-0)^{-1}F=d_{0}R_{N+1}$ for some constant $d_{0}$ (actually $F=-d_{0}N_{N+1,0}$ with $N_{N+1,0}$
defined in Corollary ~\ref{Co:admissibility2} below),
the matrix $\sigma_{1}$ is defined as $\left(
\begin{array}{lll}
0&-1\\
1&0
\end{array}
\right).$
\item[(D)]
\begin{equation}\label{eq:pkmn}
p_{k}(z,\bar{z})=\displaystyle\sum_{2N+1\geq m+n\geq
2}P^{(k)}_{m,n}(\lambda)z^{m}\bar{z}^{n}, \ k=1,2,3,4.
\end{equation} If $m,n\leq N$ then $P^{(1)}_{m,n}(\lambda)$ and $P^{(3)}_{m,n}(\lambda)$ are real,
$P^{(2)}_{m,n}(\lambda)$ and $P^{(4)}_{m,n}(\lambda)$ are purely
imaginary.
\end{enumerate}
The term $Remainder$ is bounded as
\begin{equation}\label{remainder}
\begin{array}{lll}
|Remainder(t)|&\lesssim& |z|^{2N+2}(t)+\|\langle
x\rangle^{-\mu}R_{N}(t)\|_{2}^{2}+\|R_{N}(t)\|_{\infty}^{2}\\
& &+|z|(t)\|\langle x\rangle^{-\mu}R_{2N}(t)\|_{2}
\end{array}
\end{equation} with $\nu$ being the same in Theorem
~\ref{GOLD:maintheorem}.
\end{theorem}
In the next we present the properties of the remainder function $R$ in detail.
First we define the following functions to
measure the asymptotic behavior of the parameters and remainders:
\begin{equation}\label{majorant}
\begin{array}{lll}
Z(T):=\displaystyle\max_{t\leq T}(T_{0}+t)^{\frac{1}{2N}}|z(t)|,& &
\mathcal{R}_{1}(T):=\displaystyle\max_{t\leq
T}(T_{0}+t)^{\frac{N+1}{2N}}\|\rho_{-\nu} R_{N}\|_{\mathcal{H}^{l}},\\
\mathcal{R}_{2}(T):=\displaystyle\max_{t\leq
T}(T_{0}+t)^{\frac{N+1}{2N}}\|R_{N}(t)\|_{\infty},& &
\mathcal{R}_{3}(T):=\displaystyle\max_{t\leq
T}(T_{0}+t)^{\frac{2N+1}{2N}}\|\rho_{-\nu}R_{2N}(t)\|_{2}\\
\mathcal{R}_{4}(T):=\displaystyle\max_{t\leq
T}\|R_{N}(t)\|_{\mathcal{H}^{l}}
\end{array}
\end{equation} where $l:=[\frac{d}{2}]+3,$
$T_{0}:=(|z_{1}^{(0)}|+|z_{2}^{(0)}|)^{-1}$, $\rho_{\nu}=\langle
x\rangle^{\nu}$, recall the definitions of $z_{1}^{(0)}$ and
$z_{2}^{(0)}$ and of the large constant $\nu$ in Theorem
~\ref{GOLD:maintheorem}. Before stating the main theorem we
introduce the notion of admissible functions.
\begin{definition}\label{admissible}
A vector-function $\vec{u}:\ \mathbb{R}^{d}\rightarrow
\mathbb{C}^{2}$ is admissible if the vector-function $\left(
\begin{array}{lll}
1&0\\
0&i
\end{array}
\right)\vec{u}$ has real entries.
\end{definition}
\begin{theorem}\label{THM:maintheorem2} The function $R$ in Theorem ~\ref{GOLD:maintheorem}
satisfies
the equation
\begin{equation}\label{Rorthogonal}
Im\langle R,i\phi^{\lambda}\rangle=Im\langle
R,\frac{d}{d\lambda}\phi^{\lambda}\rangle=Im\langle
R,i\eta\rangle=Im\langle R,\xi\rangle=0
\end{equation}
and the function $\vec{R}$ can be decomposed as
\begin{equation}\label{expansionR1}
\vec{R}=\displaystyle\sum_{2\leq m+n\leq
k\leq 2N}R_{mn}(\lambda)z^{m}\bar{z}^{n}+R_{k}
\end{equation} where the functions $R_{mn}:\mathbb{R}^{d}\rightarrow
\mathbb{C}^{2},\ R_{k}:\mathbb{R}^{d}\rightarrow \mathbb{R}^{2}$
have the following properties:
\begin{enumerate}
 \item[(RA)]  the function $R_{m,n}\in \mathcal{L}^{2}$ is admissible, and decays exponentially fast at $\infty$ if $\max\{m,n\}\leq N;$ and
\item[(RB)] if $\max\{m,n\}>N$ then
 the functions $R_{m,n}$ are of the form
\begin{equation}\label{eq:unusual}
\prod_{k}(L(\lambda)+ik\epsilon(\lambda)+0)^{-n_{k}}P_{c}\phi_{m,n},
\end{equation}
where the function $\phi_{m,n}$ is smooth and decays exponentially
fast at $\infty,$ $0\leq \displaystyle\sum_{k}n_{k}\leq N,$ and
$2N\geq k\geq N+1;$ note that the equation (~\ref{eq:unusual}) makes
sense in some weighted $\mathcal{L}^{2}$ space, see ~\cite{GS2,
GS3};
\end{enumerate}
the function $R_{k}$ $(2\leq k\leq 2N)$ satisfies the equation
\begin{equation}\label{equationg}
\begin{array}{lll}
\frac{d}{dt}R_{k}&=&L(\lambda)R_{k}+M_{k}(z,\bar{z})R_{k}+N_{\min\{k,N\}}(R_{\min\{k,N\}},z,\bar{z})+F_{k}(z,\bar{z}),
\end{array}
\end{equation} where
\begin{enumerate}
 \item[(1)] $F_{k}(z,\bar{z})=O(|z|^{k+1})$ is a polynomial of $z$ and
 $\bar{z}$ with $\lambda$-function-valued coefficients, and each coefficient can be written as the sum of
 functions of the form (~\ref{eq:unusual}).
 \item[(2)] $M_{k}(z,\bar{z})$ is an operator defined by
 $$M_{k}(z,\bar{z}):=\dot{\gamma}P_{c}J+\dot{\lambda}P_{c\lambda}+A_{k}(z,\bar{z}),$$ where $A_{k}(z,\bar{z})$ is a $2\times 2$ matrix bounded in the matrix norm as $$|A_{k}(z,\bar{z})|\leq c|z|^{2}e^{-\epsilon_{0}|x|}.$$
 \item[(3)] $N_{N}(R_{N},z,\bar{z})$
 satisfies the estimates
 \begin{equation}\label{remainderestimate}
 \|N_{N}(R_{N},z,\bar{z})\|_{1}+\|N_{N}(R_{N},z,\bar{z})\|_{\mathcal{H}^{2}}\leq c
 (T_{0}+t)^{-\frac{2N+3}{2N}}[Z^{2}\mathcal{R}_{1}^{2}\mathcal{R}_{2}^{4}+\mathcal{R}_{2}^{5}\mathcal{R}_{4}^{2}],
 \end{equation}
 \begin{equation}\label{remainderestiamte2}
 \begin{array}{lll}
 &
 &\|(-\triangle+1)^{2}N_{N}(R_{N},z,\bar{z})\|_{1}+\|(-\triangle+1)^{2}N_{N}(R_{N},z,\bar{z})\|_{2}\\
 &\leq&c(T_{0}+t)^{-\frac{2N+2}{2N}}(Z^{2}\mathcal{R}_{1}^{2}\mathcal{R}_{4}^{2}+Z^{2}\mathcal{R}_{1}\mathcal{R}_{2}^{2}\mathcal{R}_{4}^{3}+\mathcal{R}_{2}^{3}\mathcal{R}_{4}^{4}).
 \end{array}
 \end{equation}
\end{enumerate}
\end{theorem}
\begin{flushleft}
\textbf{Proof of Theorems ~\ref{GOLD:maintheorem} and
~\ref{THM:maintheorem2}}
\end{flushleft}
By plugging Equation (~\ref{Decom}) into Equation (~\ref{NLS}) we
have
\begin{equation}\label{Eq:R}
\begin{array}{lll}
\frac{d}{dt}\vec{R}&=&L(\lambda)\vec{R}+\dot\gamma
J\vec{R}+J\vec{N}(R,p,z) +\left(
\begin{array}{lll}
(z_{2}+p_{4})\epsilon(\lambda)\xi+\dot\gamma (z_{2}+p_{4})\eta\\
p_{1}\phi^{\lambda}-\dot\gamma [\phi^{\lambda}+(z_{1}+p_{3})\xi]
\end{array}\right)\\
&-&\left(
\begin{array}{lll}
\dot\lambda\phi_{\lambda}^{\lambda}+\partial_{t}[p_{1}\phi_{\lambda}^{\lambda}]+\dot{z}_{1}\xi+\dot\lambda[z_{1}+p_{3}]\xi_{\lambda}\\
\partial_{t}[p_{2}\phi^{\lambda}]+\dot{z}_{2}\eta+\dot\lambda
(z_{2}+p_{4})\eta_{\lambda}+(z_{1}+p_{3})\epsilon(\lambda)\eta
\end{array}
\right),
\end{array}
\end{equation}
with $\vec{R}=:\left(
\begin{array}{lll}
R_{1}\\
R_{2}
\end{array}
\right)$, $R_{1}:=Re R,$ $R_{2}:=ImR,$ $J\vec{N}(R,p,z):=\left(
\begin{array}{lll}
ImN(R,p,z)\\
-ReN(R,p,z)
\end{array}
\right)$ and
$$
\begin{array}{lll}
ImN(R,p,z):= f(|\phi^{\lambda}+I_{1}+iI_{2}|^{2})I_{2}-f((\phi^{\lambda})^{2})I_{2}
\end{array}
$$
$$
\begin{array}{lll}
ReN(R,p,z):= [f(|\phi^{\lambda}+I_{1}+iI_{2}|^{2})-f((\phi^{\lambda})^{2})](\phi^{\lambda}+I_{1})-2f^{'}((\phi^{\lambda})^{2})(\phi^{\lambda})^{2}I_{1}
\end{array}
$$ with $$I_{1}:=p_{1}\phi^{\lambda}_{\lambda}+(z_{1}+p_{3})\xi+R_{1}, \ \ I_{2}=p_{2}\phi^{\lambda}+(z_{2}+p_{4})\eta+R_{2}.$$

From Equations (~\ref{Eq:R}) and the orthogonality condition
(~\ref{Rorthogonal}) we obtain equations for $\dot{z}_{1},\
\dot{z}_{2},\ \dot\lambda,\ \dot\gamma$
\begin{equation}\label{eq:z1}
[\frac{d}{dt}({z}_{1}+p_{3})-\epsilon(\lambda)(z_{2}+p_{4})]\langle
\xi,\eta\rangle-\langle ImN(R,p,z),\eta\rangle=F_{1};
\end{equation}
\begin{equation}\label{eq:z2}
[\frac{d}{dt}({z}_{2}+p_{4})+\epsilon(\lambda)(z_{1}+p_{3})]\langle
\xi,\eta\rangle+\langle ReN(R,p,z),\xi\rangle=F_{2};
\end{equation}
\begin{equation}
\dot\gamma+\partial_{t}p_{2}-p_{1}+\frac{1}{\langle
\phi^{\lambda},\phi^{\lambda}_{\lambda}\rangle}\langle
ReN(R,p,z),\phi^{\lambda}_{\lambda}\rangle=F_{3};
\end{equation}
\begin{equation}\label{eq:lambda}
\dot\lambda+\partial_{t}p_{1}-\frac{1}{\langle
\phi^{\lambda},\phi^{\lambda}_{\lambda}\rangle}\langle
ImN(R,p,z),\phi^{\lambda}\rangle=F_{4},
\end{equation}
where the scalar functions $F_{n},\ n=1,2,3,4,$ are defined as
$$F_{1}:=\dot\gamma\langle
(z_{2}+p_{4})\eta,\eta\rangle-\dot\lambda\langle
(z_{1}+p_{3})\xi_{\lambda},\eta\rangle-\dot\gamma\langle
R_{2},\eta\rangle+\dot\lambda\langle R_{1},\eta_{\lambda}\rangle;$$
$$F_{2}:=-\dot\gamma\langle
(z_{1}+p_{3})\xi,\xi\rangle-\dot\lambda\langle
(z_{2}+p_{4})\eta_{\lambda},\xi\rangle+\dot\gamma\langle
R_{1},\xi\rangle+\dot\lambda\langle R_{2},\xi_{\lambda}\rangle;$$
$$F_{3}:=\frac{1}{\langle
\phi^{\lambda},\phi^{\lambda}_{\lambda}\rangle}[\dot\lambda \langle
R_{2},\phi_{\lambda\lambda}^{\lambda}\rangle -\dot\gamma\langle
R_{2},\phi_{\lambda}^{\lambda}\rangle-p_{2}\dot\lambda\langle\phi^{\lambda}_{\lambda},\phi^{\lambda}_{\lambda}\rangle];$$
$$F_{4}:=\frac{1}{\langle
\phi^{\lambda},\phi^{\lambda}_{\lambda}\rangle}[\dot\lambda\langle
R_{1},\phi_{\lambda}^{\lambda}\rangle +\dot\gamma\langle
R_{1},\phi^{\lambda}\rangle-p_{1}\dot\lambda\langle
\phi^{\lambda}_{\lambda\lambda},\phi^{\lambda}\rangle].$$ Moreover
using the fact that
$P_{c}\frac{d}{dt}\vec{R}=\frac{d}{dt}\vec{R}-\dot\lambda
P_{c\lambda}\vec{R}$ we rewrite Equation (~\ref{Eq:R}) as
\begin{equation}\label{RAfProj}
\frac{d}{dt}\vec{R}-L(\lambda)\vec{R}-P_{c}J\vec{N}(R,p,z)=\mathcal{G}
\end{equation}
where the function
$$\mathcal{G}:=-\dot\lambda P_{c\lambda}\vec{R}+\dot\gamma
P_{c}J\vec{R}+P_{c}\left(
\begin{array}{lll}
\dot\gamma (z_{2}+p_{4})\eta-\dot{\lambda}p_{1}\partial_{\lambda}^{2}\phi^{\lambda}-\dot\lambda (z_{1}+p_{3})\xi_{\lambda}\\
-\dot{\lambda}p_{2}\phi^{\lambda}_{\lambda}-\dot\lambda
(z_{2}+p_{4})\eta_{\lambda}
\end{array}
\right).$$

In the next we sketch the proof of the theorem except Equation
(~\ref{FGR}). We start with proving the expansion of $R$ and the
expressions for $\dot{z},\ \dot\lambda,\ \dot\gamma$ and $p_{l},\
l=1,2,3,4,$ by induction on $k=m+n$, $k=2,3,\cdot\cdot\cdot,$ with
the main tool Lemma ~\ref{LM:admissibility1} below. For the detail
of finding $R_{m,n},\ \Lambda_{m,n},\ \Gamma_{m,n},\ Z_{m,n}$ we
refer to ~\cite{GS2, GS3}. For $k=2,$ by Lemma
~\ref{LM:admissibility1} we have that $iN_{m,n}^{(2)},$ defined in
(~\ref{NonExpan}) below, is admissible (Note $N_{m,n}^{(2)},\
m+n=2,$ does not depend on any $R_{m',n'}, m+n\geq 2$ by Lemma
~\ref{LM:admissibility1} and the term $\mathcal{G}$ has no
contribution for terms of order $|z|^{2}$). Thus
$R_{m,n}=-[L(\lambda)+i(m-n)\epsilon(\lambda)]^{-1}N_{m,n}^{(2)}$,
$m+n=2,$ is admissible if $N>1$. Moreover by the properties of
$Z_{m,n},\ \Lambda_{m,n}$ and $\Gamma_{m,n}$ in the equations for
$\dot\lambda,\ \dot\gamma,\ \dot{z}=\dot{z}_{1}+i\dot{z}_{2}$ in
(~\ref{ExpanLambda})-(~\ref{FGRZ}), the equations
(~\ref{eq:z1})-(~\ref{eq:lambda}) and the observation that the term
$F_{j},\ j=1,2,3,4$ has no terms of order $|z|^{2}$ we have if
$m+n=2$ and $m\not=n+1$ then
$$
\begin{array}{lll}
& &[-i(m-n)\epsilon(\lambda)P_{m,n}^{(3)}-\epsilon(\lambda)P_{m,n}^{(4)}]+i[-i(m-n)\epsilon(\lambda)P_{m,n}^{(4)}+\epsilon(\lambda)P_{m,n}^{(3)}]\\
&=&\frac{1}{\langle \xi,\eta\rangle}[\langle
ImN_{m,n}^{(2)},\eta\rangle-i\langle ReN_{m,n}^{(2)},\xi\rangle]
\end{array}$$
and if $m'+n'=2$ and $m'\not=n'$ then
$$-i(m'-n')\epsilon(\lambda)P^{(2)}_{m',n'}-P^{(1)}_{m',n'}=-\frac{1}{\langle \phi^{\lambda},\phi^{\lambda}_{\lambda}\rangle}\langle ReN_{m',n'}^{(2)},\phi^{\lambda}_{\lambda}\rangle,$$
$$-i(m'-n')\epsilon(\lambda)P^{(1)}_{m',n'}=\frac{1}{\langle \phi^{\lambda},\phi^{\lambda}_{\lambda}\rangle}\langle ImN^{(2)}_{m',n'},\phi^{\lambda}\rangle$$
where $\left(
\begin{array}{lll}
ImN^{(2)}_{m,n}\\
-ReN^{(2)}_{m,n}
\end{array}
\right):=N_{m,n}^{(2)}.$ We observe the solutions satisfying the
condition $P^{(k)}_{m,n}=\overline{P^{(k)}_{n,m}}$, $k=1,2,3,4,$ to
the equations above exist, and by the properties of $N_{m,n}^{(2)}$
we obtain $P^{(1)}_{m,n}(\lambda)$ and $P^{(3)}_{m,n}(\lambda)$ are
real, $P^{(2)}_{m,n}(\lambda)$ and $P^{(4)}_{m,n}(\lambda)$ are
purely imaginary. Thus we finish the proof of the step of $m+n=2.$
For the step $m+n=3,$ $P^{(l)}_{m',n'}$, $R_{m',n'},\ m'+n'=2$ and
Lemma ~\ref{LM:admissibility1} enable us to determine
$P^{(l)}_{m,n},\ R_{m,n},\ m+n=3$ and prove all the properties in
the theorem. By using the procedures above repeatedly we finish the
proof.

By using the same techniques in ~\cite{GS2,GS3} and Equation
(~\ref{equationg}) we prove the decay of $R_{N}$ and $R_{2N}$ in
appropriate norms; and by (~\ref{ExpanLambda}) there exists a new
parameter $\beta=z+O(|z|^{2})$ such that
\begin{equation}\label{FGR2}
\dot\beta=-i\epsilon(\lambda)\beta+\sum_{1\leq n\leq
N}Y_{n}\beta^{n+1}\bar\beta^{n}+Remainder,
\end{equation}
where $Y_{n}$ is purely imaginary if $n<N$ and $Re Y_{N}=Re
Z_{N+1,N}$, consequently we have
$$\frac{1}{2}\frac{d}{dt}|z|^{2}=ReZ_{N+1,N}|z|^{2(N+1)}+|z|Remainder$$ which implies the decay of $\beta$, hence $y$,
by the assumption $ReZ_{N+1,N}<0$ and the estimates on
$|z|Remainder;$ and similar as in ~\cite{GS3} by the expressions for
$\dot\lambda$ in (~\ref{ExpanLambda}) there exist $c_{m,n}(\lambda)$
such that
$$\frac{d}{dt}[\lambda-\sum_{2N+1\geq m+n\geq
2}c_{m,n}(\lambda)z^{m}\bar{z}^{n}]=Remainder$$ with the the term
$Remainder$ satisfies (~\ref{remainder}), consequently $\lambda$ is
convergent due to the decay of $z$ and the fact that $Remainder(t)$
is integrable at $\infty.$

The proof of Equation (~\ref{FGR}) is in Sections
~\ref{SEC:ProofFGR} and ~\ref{sec:CaseN3} below.
\begin{flushright}
$\square$
\end{flushright}

In the proof we use the following lemma.
\begin{lemma}\label{LM:admissibility1}
\begin{itemize}
\item[(A)]
If $G_{1}$ is a vector-function from $\mathbb{R}^{d}$ to $
\mathbb{C}^{2}$ such that $iG_{1}$ is admissible and $G_{1}\in
\mathcal{L}^{2}$, then the vector function
$(L(\lambda)+i\mu)^{-1}P_{c}G_{1}$ is admissible for any $\mu\in
(-\lambda,\lambda).$
\item[(B)]
Let $R_{k},\ R_{m,n},\ p_{v}$ be the same as those in
(~\ref{expansionR1}) and (~\ref{Decom}) respectively and
$$p_{v}=\displaystyle\sum_{m+n\geq
2}q_{m,n}^{(v)}(\lambda)z^{m}\bar{z}^{n}, \ v=1,2,3,4$$ with
$q^{(1)}_{m,n}(\lambda)$ and $q^{(3)}_{m,n}(\lambda)$ are real,
$q^{(2)}_{m,n}(\lambda)$ and $q^{(4)}_{m,n}(\lambda)$ are purely
imaginary if $m,\ n\leq \min\{k,N\}$. Then $J \vec{N}(R,p,z)$ can be
expanded as
\begin{equation}\label{NonExpan}
\begin{array}{lll}
J\vec{N}(R,p,z)&=&\displaystyle\sum_{2\leq m+n\leq
2N}z^{m}\bar{z}^{n}N^{(k)}_{mn}(\lambda)+A_{k}(z,\bar{z})R_{k}\\
& &+N_{\min\{k,N\}}(R_{\min\{k,N\}},p,z)+Remainder1
\end{array}
\end{equation}
where $iN_{m,n}^{(k)}(\lambda):\mathbb{R}^{3}\rightarrow
\mathbb{C}^{2}$ is admissible if $m,n\leq N$ and only depends on
$R_{m',n'}$ $q^{(v)}_{m',n'}$ with $m'\leq m, n'\leq n$ and
$|m'-m|+|n'-n|\not=0$; the $2\times 2$-matrix function
$A_{k}(z,\bar{z})$ is bounded in the matrix norm as
$|A_{k}(z,\bar{z})|\leq c|z||e^{-\epsilon_{0}|x|};$
$N_{j}(R_{j},p,z)$, $j\leq N,$ contains all the nonlinear terms in
$R_{j}$, $N_{N}(R_{N},p,z)$ admits the same estimates as in
(~\ref{remainderestimate}) (~\ref{remainderestiamte2}); and the term
$Remainder_{1}$ has the estimate
\begin{equation}\label{remainder1}
|Remainder_{1}|\leq c|z|^{2N+1}e^{-\epsilon_{0}|x|}.
\end{equation}
\end{itemize}
\end{lemma}
\begin{proof}
The first part is copied from ~\cite{GS2, GS3}.

For the second part we prove the expansion for
$J\vec{N}(\vec{R},p=0,z)$ in ~\cite{GS2, GS3}. The proof is almost
identical to that, thus omitted.
\end{proof}
The following corollary will be used later, whose proof is obvious
by (~\ref{NonExpan}), thus omitted.
\begin{corollary}\label{Co:admissibility2} Let $k=2N$ in (~\ref{expansionR1}), then
$$J
\vec{N}(R,p,z)=\displaystyle\sum_{m+n=2}^{2N}N_{m,n}(\lambda)z^{m}\bar{z}^{n}+N_{N}(R_{N},p,z)+Remainder1$$
with $iN_{m,n}$ admissible for $m,n\leq N$; $N_{N}(R_{N},p,z)$ and
$Remainder1$ have the same estimates as above.
\end{corollary}
\section{The Effective Equations for $\dot{z},\ \dot\lambda,\
\dot\gamma$ and $ R$}\label{SEC:effective} The equations for the
parameters $\dot{z},\ \dot\lambda,\ \dot\gamma$ and the function $R$
have some terms having no role on determining $ReZ_{N+1,N}$. In this
subsection we separate the effective parts from the "useless" parts.
We will use the following lemma to show that the terms on the right
hand side of Equations (~\ref{eq:z1})-(~\ref{eq:lambda}) and
(~\ref{RAfProj}) are useless.
\begin{lemma}\label{LM:junk}
If we expand the scalar functions $F_{1}+iF_{2},$ $F_{3},$ $F_{4}$
and the function $\mathcal{G}$ in $z$ and $\bar{z},$ then we have
\begin{equation}\label{junk}
F_{1}+iF_{2}=\sum_{2\leq m+n\leq
2N+1}K^{(1)}_{m,n}(\lambda)z^{m}\bar{z}^{n}+Remainder
\end{equation}
$$F_{3}=\sum_{2\leq m+n\leq 2N}K^{(2)}_{m,n}(\lambda)z^{m}\bar{z}^{n}+Remainder$$
$$F_{4}=\sum_{2\leq m+n\leq 2N}K^{(3)}_{m,n}(\lambda)z^{m}\bar{z}^{n}+Remainder$$
\begin{equation}\label{junkG}
\mathcal{G}=\sum_{2\leq m+n\leq
2N}\mathcal{G}_{m,n}(\lambda)z^{m}\bar{z}^{n}+GN(R_{N},z,\bar{z})+Remainder1
\end{equation} where
if $m,n\leq N+1$, then the coefficients $K^{(1)}_{m,n},\
K^{(3)}_{m,n}$ are purely imaginary, $K^{(2)}_{m,n}$ is real, the
function $i\mathcal{G}_{m,n}$ is admissible; $K^{(k)}_{m,n}=0$ and
$\mathcal{G}_{m,n}=0$ if $(m,n)=(0,l)\ \text{or}\ (l,0)$ with $l\leq
2N;$ moveover the function $Remainder1$ satisfies the estimate
(~\ref{remainder1}), the term $GN(R_{N},z,\bar{z})$ admits the same
estimates as $N_{N}(R_{N},z,\bar{z})$ in (~\ref{remainderestimate})
(~\ref{remainderestiamte2}).
\end{lemma}
\begin{proof}
By the expansions for the scalar functions $\dot\lambda,$
$\dot\gamma$ and the function $\vec{R}$ in Theorems
~\ref{GOLD:maintheorem} and ~\ref{THM:maintheorem2}, we prove the
lemma by direct computations.
\end{proof}
\section{Proof of Equation (~\ref{FGR}) when $N=2$}\label{SEC:ProofFGR}
In this section we prove Equation (~\ref{FGR}) when $N=2.$ For
consideration of space we only prove the case of the nonlinearity
$f(x)=x.$ One can show that the results hold for more general
nonlinearities, but the computations are much messier. The main
result is
\begin{theorem}\label{THM:mainresultN2}
If
$N=2$ and the nonlinearity $f(x)=x$ in (~\ref{NLS}), then
\begin{equation}\label{eq:ComZ32}
ReZ_{3,2}=\frac{6}{\langle \xi,\eta\rangle}Im\langle
\sigma_{1}(L(\lambda)-3i\epsilon(\lambda)-0)^{-1}P_{c}N_{3,0},N_{3,0}\rangle,
\end{equation} where, recall the definition of $N_{3,0}$ from
Corollary ~\ref{Co:admissibility2}.
\end{theorem}
The proof is in subsection ~\ref{subsec:proofN2}.
\subsection{the coefficient of $z^{3}\bar{z}^{2}$ in the expansion of $\langle Im
N(R,p,z),\eta\rangle-i\langle ReN(R,p,z),\xi\rangle$} Suppose that
the expansion of $\langle Im N(R,p,z),\eta\rangle-i\langle
ReN(R,p,z),\xi\rangle$, the part on the left hand side of Equations
(~\ref{eq:z1}) and (~\ref{eq:z2}), in $z$ and $\bar{z}$ is
$$\langle
ImN(R,p,z),\eta\rangle-i\langle ReN(R,p,z),\xi\rangle=\sum_{2\leq
m+n\leq 5}X_{m,n}z^{m}\bar{z}^{n}+Remainder.$$ In the next we
compute $X_{3,2}$. Before starting the proof we define the following
functions and constants:
$$K_{1}:=\left(
\begin{array}{lll}
-P_{2,0}^{(2)}(\phi^{\lambda})^{2}\xi+iP_{2,0}^{(1)}\phi^{\lambda}\phi^{\lambda}_{\lambda}\eta\\
3
P_{2,0}^{(1)}\phi^{\lambda}\phi^{\lambda}_{\lambda}\xi-iP_{2,0}^{(2)}(\phi^{\lambda})^{2}\eta
\end{array}
\right),\ K_{2}:=\left(
\begin{array}{lll}
i\phi^{\lambda}\eta, &-\phi^{\lambda}\xi \\
3\phi^{\lambda}\xi, &-i\phi^{\lambda}\eta
\end{array}
\right)R_{2,0},$$
$$K_{3}:=\frac{1}{8}\left(
\begin{array}{lll}
-i\eta^{3}+i\xi^{2}\eta\\
-\xi\eta^{2}+\xi^{3}
\end{array}
\right),\ K_{4}:=\left(
\begin{array}{lll}
-\phi^{\lambda}\xi\eta[P_{2,0}^{(4)}-iP^{(3)}_{2,0}]\\
-i\phi^{\lambda}\eta^{2}P_{2,0}^{(4)}+3\phi^{\lambda}\xi^{2}P_{2,0}^{(3)}
\end{array}
\right)$$ and
$$
D_{1}(\lambda):=i[P_{3,1}^{(1)}\langle
\phi^{\lambda}(\eta^{2}-3\xi^{2}),\phi_{\lambda}^{\lambda}\rangle-iP_{3,1}^{(2)}\langle
2\xi\eta\phi^{\lambda},\phi^{\lambda}\rangle],
$$
$$D_{2}(\lambda):=-2i\langle \sigma_{1} R_{3,0},K_{1}+K_{2}+3K_{3}+K_{4}\rangle,$$
$$D_{3}(\lambda):=\langle R_{3,1},\left(
\begin{array}{lll}
-i(\phi^{\lambda}\eta^{2}-3\phi^{\lambda}\xi^{2})\\
2\phi^{\lambda}\xi\eta
\end{array}
\right)\rangle,$$ and $$D_{4}:=P_{3,1}^{(3)}[-3i\langle
\phi^{\lambda}\xi^{2},\xi\rangle+i\langle
\phi\xi\eta,\eta\rangle]+P_{3,1}^{(4)}\langle
2\phi^{\lambda}\eta^{2},\xi\rangle$$
 and the matrix $\sigma_{1}:=\left(
\begin{array}{lll}
0&-1\\
1&0
\end{array}
\right),$ recall the definitions of $P_{m,n}^{(l)},\ l=1,2,3,4$ in
Theorems ~\ref{GOLD:maintheorem}. The main result is
\begin{proposition}\label{Pro:X32}
$$ReX_{3,2}=\sum_{n=1}^{4}Re D_{n}(\lambda)=6Im\langle \sigma_{1}R_{3,0},\sum_{n=1}^{4}K_{n}\rangle.$$
\end{proposition} The proof is at the end of this section.
In the proof the following observation will be used
\begin{lemma}\label{LM:Cal}
$P_{m,n}^{(k)}=\overline{P_{n,m}^{(k)}}$ for any $m,n$ and $k$; If
$k=1,3$ and $m,n\leq 2$ then $P^{(k)}_{m,n}$ is real; if $k=2,4$ and
$m,n\leq 2$ then $P^{(k)}_{m,n}$ is purely imaginary. The constants
$P_{m,n}^{(k)}$ have the follow expressions
$$P_{2,0}^{(1)}=\frac{1}{4\epsilon(\lambda)\langle
\phi^{\lambda},\phi^{\lambda}_{\lambda}\rangle}\langle
\phi^{\lambda}\xi\eta,\phi^{\lambda}\rangle ,\
P_{2,0}^{(2)}=-\frac{P_{2,0}^{(1)}}{i2\epsilon(\lambda)}+\frac{\langle
3\phi^{\lambda}\xi^{2}-\phi^{\lambda}\eta^{2},\phi^{\lambda}_{\lambda}\rangle}{8i\epsilon(\lambda)\langle
\phi^{\lambda},\phi^{\lambda}_{\lambda}\rangle},$$
$$
-2iP_{2,0}^{(3)}-P_{2,0}^{(4)}=\frac{1}{\epsilon(\lambda)\langle
\xi,\eta\rangle}\langle
-\frac{i}{2}\phi^{\lambda}\xi\eta,\eta\rangle,
$$
$$
-2iP_{2,0}^{(4)}+P_{2,0}^{(3)}=-\frac{1}{4\epsilon(\lambda)\langle
\xi,\eta\rangle}\langle \phi^{\lambda}\xi,3\xi^{2}-\eta^{2}\rangle,
$$
$$-2iP_{3,1}^{(3)}-P_{3,1}^{(4)}=\frac{1}{\epsilon(\lambda)\langle \xi,\eta\rangle}\langle R_{3,0},\left(
\begin{array}{lll}
-i\phi^{\lambda}\eta^{2}\\
\phi^{\lambda}\xi\eta
\end{array}
\right)+iU,$$
$$-2iP_{3,1}^{(4)}+P_{3,1}^{(3)}=\frac{1}{\epsilon(\lambda)\langle \xi,\eta\rangle}\langle R_{3,0},\left(
\begin{array}{lll}
-3\phi^{\lambda}\xi^{2}\\
i\phi^{\lambda}\xi\eta
\end{array}
\right)+U,$$
$$P_{3,1}^{(1)}=\frac{1}{2\epsilon(\lambda)\delta^{'}(\lambda)}\langle
R_{3,0},\left(
\begin{array}{lll}
-(\phi^{\lambda})^{2}\eta\\
-i(\phi^{\lambda})^{2}\xi
\end{array}
\right)\rangle+U$$
$$P_{3,1}^{(2)}=-\frac{P_{3,1}^{(1)}}{i2\epsilon(\lambda)}
-\frac{i}{2\epsilon(\lambda)\delta^{'}(\lambda)}\langle
R_{3,0},\left(
\begin{array}{lll}
3\phi^{\lambda}\phi^{\lambda}_{\lambda}\xi\\
-i\phi^{\lambda}\phi_{\lambda}^{\lambda}\eta
\end{array}
\right)\rangle+iU,$$ $U$ stands for real (different) constants
throughout the rest of the paper. If $m,n\leq 3$ and
$|m-3|+|n-3|\not=0$ then
\begin{equation}\label{eq:RMN}
R_{m,n}=-(L(\lambda)+i(m-n)\epsilon(\lambda)+0)^{-1}[N_{m,n}-\mathcal{G}_{m,n}].
\end{equation}
\end{lemma}
\begin{proof} The property $P_{m,n}^{(k)}=\overline{P_{n,m}^{(k)}}$
following from the fact that $\lambda, \ \gamma,\ z_{1},\ z_{2}$ are
real functions. The fact $P^{(k)}_{m,n}$ being real or purely
imaginary is based on the properties of $N_{m,n}$ in Corollary
~\ref{Co:admissibility2}. By Equations
(~\ref{eq:z1})-(~\ref{eq:lambda}) and the properties of
$\dot\lambda,\dot\gamma,\dot{z}$ in Theorem ~\ref{GOLD:maintheorem}
we compute to get the constants $P_{m,n}^{(k)}.$

The formula of $R_{m,n}$ follows from the expansion, we refer to our
paper ~\cite{GS2, GS3}.
\end{proof}
We have the important remark for $-0$ symbol in the expression of
$R_{m,n}$.
\begin{remark}
If $|m-n|\leq 2$ then
$$(L(\lambda)+i(m-n)\epsilon(\lambda)+0)^{-1}=(L(\lambda)+i(m-n)\epsilon(\lambda))^{-1}.$$
\end{remark}
\textbf{Proof of Proposition ~\ref{Pro:X32}} For the fact
$ReX_{3,2}=\sum_{n=1}^{4}Re D_{n}(\lambda)$ we only use Lemma
~\ref{LM:Cal}, the fact that $R_{m,n}$ is admissible if $m,n\leq
N=2$ in Theorem ~\ref{THM:maintheorem2} and direct computation, thus
omit the detail.

For the fact on the right hand side we transform the expression for
$D_{n},\ n=1,2,3,4.$

First we prove the $D_{1}$ part. By the expressions of
$P_{2,0}^{(1)},\ P_{3,1}^{(1)},\ P_{2,0}^{(2)},\ P_{3,1}^{(2)}$ and
the properties of $P_{m,n}^{(k)}$ in Lemma ~\ref{LM:Cal} we have
$$
\begin{array}{lll}
ReD_{1}&=&Rei[P_{3,1}^{(1)}\langle
\phi^{\lambda}(\eta^{2}-3\xi^{2}),\phi_{\lambda}^{\lambda}\rangle+2P_{3,1}^{(2)}\langle
-i\xi\eta\phi^{\lambda},\phi^{\lambda}\rangle] \\
&=&4Re\langle R_{3,0}, \left(
\begin{array}{lll}
3P_{2,0}^{(1)}i\phi^{\lambda}\phi^{\lambda}_{\lambda}\xi+P_{2,0}^{(2)}(\phi^{\lambda})^{2}\eta\\
P_{2,0}^{(1)}\phi^{\lambda}\phi^{\lambda}_{\lambda}\eta+i
P_{2,0}^{(2)}(\phi^{\lambda})^{2}\xi
\end{array}
\right)\rangle\\
&=&-Re\langle R_{3,0}, 4i\sigma_{1} K_{1}\rangle.
\end{array}
$$

For $D_{3}$ we compute $R_{2,0}$ and $R_{3,1}$ first. Recall the
forms of $R_{m,n}$ in Lemma ~\ref{LM:Cal} and the properties of
$\mathcal{G}_{m,n}$ in (~\ref{junkG}). First we have
\begin{equation}\label{eq:r20}
R_{2,0}=-[L(\lambda)+2i\epsilon(\lambda)]^{-1}P_{c}N_{2,0}
\end{equation} with $N_{2,0}=\frac{1}{4}\left(
\begin{array}{lll}
-2i\phi^{\lambda}\xi\eta\\
-3\phi^{\lambda}\xi^{2}+\phi^{\lambda}\eta^{2}
\end{array}
\right).$ The function $R_{3,1}$ has two parts
\begin{equation}\label{R31}
R_{3,1}=-[L(\lambda)+2i\epsilon(\lambda)]^{-1}P_{c}(G_{1}+K_{5})
\end{equation} where the vector function $iG_{1}$
is admissible, hence
$-[L(\lambda)+2i\epsilon(\lambda)]^{-1}P_{c}G_{1}$ is admissible
(thus "useless") by Lemma ~\ref{LM:admissibility1} and the vector
function $K_{5}$ is
\begin{equation}\label{eq:k4}
K_{5}:=-\left(
\begin{array}{lll}
i\phi^{\lambda}\eta & \phi^{\lambda}\xi \\
-3\phi^{\lambda}\xi &-i\phi^{\lambda}\eta
\end{array}
\right)R_{3,0}
\end{equation}
which is part of $N_{3,1}.$ By the observations that
$\sigma_{1}(L(\lambda))^{*}\sigma_{1}=L(\lambda),$
$\sigma^{2}_{1}=-1$ and the definition of $K_{2}$ we have
$$
ReD_{3} =-Re\langle R_{3,0},4i\sigma_{1}K_{2}\rangle.
$$

$D_{4}$ admits the form $ReD_{4}(\lambda)=4Im\langle
\sigma_{1}R_{3,0},K_{4}\rangle,$ which follows from some
manipulation on the expressions of $P_{2,0}^{(n)}$ and
$P_{3,1}^{(n)}$. The proof is tedious, but not hard, thus omitted.

Collecting the computation above and recalling the form of $D_{2}$
we have $$\sum_{n=1}^{4}ReD_{n}(\lambda)=6Im\langle
\sigma_{1}R_{3,0},\sum_{n=1}^{4}K_{n}\rangle.$$ Thus the proof is
complete.
\begin{flushright}
$\square$
\end{flushright}
\subsection{Proof of Theorem ~\ref{THM:mainresultN2}}\label{subsec:proofN2}
\begin{proof}
Recall that $X_{3,2}$ is the coefficient of $z^{3}\bar{z}^{2}$ in
the expansion of the term $\langle ImN(R,p,z),\eta\rangle-i\langle
ReN(R,p,z),\xi\rangle.$ First by Equations (~\ref{eq:z1}),
(~\ref{eq:z2}) and (~\ref{junk}) we have $ReZ_{3,2}=\frac{1}{\langle
\xi,\eta\rangle}ReX_{3,2}$, then by Proposition ~\ref{Pro:X32}
\begin{equation}\label{eq:Claim3D}
ReX_{3,2}=\sum_{n=1}^{4}ReD_{n}(\lambda)=6Im\langle
\sigma_{1}R_{3,0},\sum_{n=1}^{4}K_{n}\rangle,
\end{equation} where, recall the definitions of $K_{n}$ before (~\ref{eq:ComZ32}).

By (~\ref{eq:RMN}) and the properties of $\mathcal{G}$ in
(~\ref{junkG}) for $N=2$
$$
R_{3,0}=-(L(\lambda)+3i\epsilon(\lambda)+0)^{-1}P_{c}N_{3,0} $$ and
we compute to get
$$N_{3,0}=-\sum_{n=1}^{4}K_{n}.$$
Consequently we have
$$
\begin{array}{lll}
ReZ_{3,2}&=&\frac{1}{\langle
\xi,\eta\rangle}\displaystyle\sum_{n=1}^{4}ReD_{n}(\lambda)\\
&=&\frac{6}{\langle \xi,\eta\rangle}Im\langle
\sigma_{1}(L(\lambda)+3i\epsilon(\lambda)+0)^{-1}P_{c}N_{3,0},N_{3,0}\rangle\\
&=&\frac{6}{\langle \xi,\eta\rangle}Im\langle
\sigma_{1}(L(\lambda)-3i\epsilon(\lambda)-0)^{-1}P_{c}N_{3,0},N_{3,0}\rangle.
\end{array}
$$ Thus the proof is complete.
\end{proof}
\section{Proof of Equation (~\ref{FGR}) when
$N=3$}\label{sec:CaseN3}
The proof of $N=3$ is more involved, but the ideas are the same.
Thus in this section we just state the key steps and only prove the case $f(x)=x$ and the potential $V$ is spherical symmetric. The following is
the main result. Recall that $\sigma_{1}=\left(
\begin{array}{lll}
0&-1\\
1&0
\end{array}
\right)$ and the definition of $N_{4,0}$ from Corollary
~\ref{Co:admissibility2}.
\begin{theorem}\label{THM:mainresult}
If $N=3$ and the nonlinearity $f(x)=x$ and the potential $V$ is spherical symmetric in (~\ref{NLS}), then the
constant $Z_{4,3}$ satisfies the equation
$$ReZ_{4,3}:=
\frac{8}{\langle \xi,\eta\rangle}Im\langle
\sigma_{1}(L(\lambda)-4i\epsilon(\lambda)-0)^{-1}N_{4,0},N_{4,0}\rangle$$
where, recall the definition of $N_{4,0}$ from Corollary
~\ref{Co:admissibility2} (for N=3).
\end{theorem}
In the next we state the following observations which will simplify
our computation. We defined two spaces $$\mathbb{S}=\{g:\ \mathbb{R}^{d}\rightarrow
\mathbb{C}|\ g(-x)=g(x)\};\ \ \mathbb{AS}=\{g:\
\mathbb{R}^{d}\rightarrow \mathbb{C}|\ g(-x)=-g(x)\}.$$ Obviously
$S\perp AS$ and $\phi^{\lambda},\phi^{\lambda}_{\lambda}\in S$ and
$\xi,\eta\in AS$ (see Corollary ~\ref{mainpo1}). Recall that $N_{m,n}$ and $\mathcal{G}_{m,n}$ are
defined in Corollary ~\ref{Co:admissibility2} and Equation
(~\ref{junkG}) respectively.
\begin{lemma}\label{LM:keyfact} If the potential $V$ is spherical symmetric, then
$R_{m,n},\ N_{m,n}\in \mathbb{S}$ if $m+n=2,4,6;$ $R_{m,n},\
N_{m,n}\in \mathbb{AS}$ if $m+n=3,5.$ Moreover
\begin{equation}\label{eq:genefacts}
P^{(3)}_{m,n}=P^{(4)}_{m,n}=P^{(1)}_{m',n'}=P^{(2)}_{m',n'}=0
\end{equation} if
$m+n$ is even or $m'+n'$ is odd. If $m,n\leq 3$ then $P^{(1)}_{m,n},
\ P^{(3)}_{m,n}$ are real and $P^{(2)}_{m,n}, \ P^{(4)}_{m,n}$ are
purely imaginary, and
\begin{equation}\label{eq:conj}
P^{(k)}_{m,n}=P^{(k)}_{n,m}\ \text{if}\ k=1,3,\
P^{(l)}_{m,n}=-P^{(l)}_{n,m}\ \text{if}\ l=2,4.
\end{equation}

If $m,n\leq 3$ and $|m-3|+|n-3|\not=0$ then
$$R_{m,n}=-(L(\lambda)+i(m-n)\epsilon(\lambda)+0)^{-1}[N_{m,n}-\mathcal{G}_{m,n}].$$
\end{lemma}
\begin{proof} The facts $P_{m,n}^{(k)}$ being real or purely imaginary follow from similar discussion in Lemma ~\ref{LM:Cal}.
The rest follows from direct computation based on the properties of $N_{m,n}$ in Corollary ~\ref{Co:admissibility2} and the observations made on $\xi,\eta,\phi^{\lambda},\phi^{\lambda}_{\lambda}$ before the lemma, thus we omit the detail.
\end{proof}
\subsection{The Explicit Forms of $R_{4,0}$ and $N_{4,0}$}  We denote the coefficients of $z^{m}\bar{z}^{n}$ in
the expansion of $ImN(R,p,z),\ (ReN(R,p,z))$ as $ImN_{m,n},\
(ReN_{m,n})$. The following lemma is the main result of this
section.
\begin{lemma}\label{LM:R40}
\begin{equation}\label{eq:defineR40}
R_{4,0}=-(L(\lambda)+4i\epsilon(\lambda)+0)^{-1}N_{4,0}
\end{equation} with
$N_{4,0}$ given explicitly below.
\end{lemma}
The proof follows from the definition of $R_{4,0}$ in Lemma
~\ref{LM:keyfact}, the property $\mathcal{G}_{4,0}=0$ in
(~\ref{junkG}) for $N=3$.

In the next we compute $N_{4,0}=\left(
\begin{array}{lll}
ImN_{4,0}\\
-ReN_{4,0}
\end{array}
\right)$ explicitly. As usual we denote $R_{m,n}=\left(
\begin{array}{lll}
R_{m,n}^{(1)}\\
R_{m,n}^{(2)}
\end{array}
\right)$ and recall that for some $m,n,k$, $P_{m,n}^{(k)}=0$ from
Lemma ~\ref{LM:keyfact}. We compute to get
\begin{equation}\label{eq:ImN40}
ImN_{4,0}= \displaystyle\sum_{n=1}^{3}\mathcal{H}_{n}^{(1)},\
ReN_{4,0}= \displaystyle\sum_{n=1}^{3} \mathcal{H}^{(2)}_{n}
\end{equation}
with the functions $\mathcal{H}_{n}^{(k)},\ n=1,2,3,$ $k=1,2,$
defined as
$$
\begin{array}{lll}
\mathcal{H}_{1}^{(1)}&:=&2(\phi^{\lambda})^{2}R_{2,0}^{(1)}P_{2,0}^{(2)}+2\phi^{\lambda}R_{2,0}^{(1)}R_{2,0}^{(2)}-\frac{i}{2}\xi\eta
R_{2,0}^{(1)}
\\
&
&+2\phi^{\lambda}\phi^{\lambda}_{\lambda}P_{2,0}^{(1)}R_{2,0}^{(2)}+\frac{1}{4}\xi^{2}R_{2,0}^{(2)}-\frac{3}{4}\eta^{2}R_{2,0}^{(2)}\\
&
&+2(\phi^{\lambda})^{2}\phi^{\lambda}_{\lambda}P_{2,0}^{(1)}P_{2,0}^{(2)}+\frac{1}{4}\phi^{\lambda}\xi^{2}P_{2,0}^{(2)}
-\frac{3}{4}\phi^{\lambda}\eta^{2}P_{2,0}^{(2)}-\frac{i}{2}\phi^{\lambda}_{\lambda}\xi\eta
P_{2,0}^{(1)};
\end{array}
$$
$$\mathcal{H}_{2}^{(1)}:=\phi^{\lambda}\xi\eta
[P_{3,0}^{(4)}-iP^{(3)}_{3,0}];$$
$$\mathcal{H}_{3}^{(1)}:=\phi^{\lambda}\xi R_{3,0}^{(2)}-i\phi^{\lambda}\eta R_{3,0}^{(1)};$$
$$
\begin{array}{lll}
\mathcal{H}_{1}^{(2)} &=
&3\phi^{\lambda}R_{2,0}^{(1)}R_{2,0}^{(1)}+6\phi^{\lambda}\phi^{\lambda}_{\lambda}P_{2,0}^{(1)}R_{2,0}^{(1)}
+\frac{3}{4}\xi^{2}R_{2,0}^{(1)}-\frac{1}{4}\eta^{2}R_{2,0}^{(1)}\\&
& +\phi^{\lambda}R_{2,0}^{(2)}R_{2,0}^{(2)}
+2(\phi^{\lambda})^{2}R_{2,0}^{(2)}P_{2,0}^{(2)}-\frac{i}{2}\xi\eta
R_{2,0}^{(2)}\\
&
&+(\phi^{\lambda})^{3}P_{2,0}^{(2)}P_{2,0}^{(2)}+3\phi^{\lambda}(\phi^{\lambda}_{\lambda})^{2}P_{2,0}^{(1)}P_{2,0}^{(1)}\\
& &-\frac{i}{2}\phi^{\lambda}\xi\eta
P_{2,0}^{(2)}+\frac{3}{4}\phi^{\lambda}_{\lambda}\xi^{2}P_{2,0}^{(1)}-\frac{1}{4}\phi^{\lambda}_{\lambda}\eta^{2}P_{2,0}^{(1)};
\end{array}
$$
$$\mathcal{H}_{2}^{(2)}:=3\phi^{\lambda}\xi^{2}P_{3,0}^{(3)}-i\phi^{\lambda}\eta^{2}P_{3,0}^{(4)};$$
$$\mathcal{H}_{3}^{(2)}:=3\phi^{\lambda}\xi R_{3,0}^{(1)}-i\phi^{\lambda}\eta R_{3,0}^{(2)}.$$

\subsection{The Explicit Form of $Z_{4,3}$ in (~\ref{FGRZ})}
The following lemma is the main result of this section.
\begin{lemma}\label{LM:z43}
\begin{equation}\label{eq:expZ43}
ReZ_{4,3}=-\frac{1}{\langle \xi,\eta\rangle}Re[E_{1}+E_{2}+E_{3}]
\end{equation}
where the constants $E_{n},\ n=1,2,3,$ are defined below.
\end{lemma}
\begin{proof}
By Lemma ~\ref{LM:junk} we have that the right hand sides of
(~\ref{eq:z1}) and (~\ref{eq:z2}) are useless in determining the
real part of $Z_{4,3}.$ Thus $ReZ_{4,3}=\frac{1}{\langle
\xi,\eta\rangle}Re[\langle ImN_{4,3},\eta\rangle-i\langle
ReN_{4,3},\xi\rangle]$ by the fact $z=z_{1}+i z_{2}.$ What is left
is just direct computation while deserting useless terms by using
the property of $P_{m,n}^{(k)},\ k=1,2,3,4$ in Lemma
~\ref{LM:keyfact} and using the fact $R_{m,n}, \ m,n\leq 3$ are
admissible.
\end{proof}
In the next we give the definitions of $E_{n}, \ n=1,2,3,$ and the
various related terms.
$$E_{1}:=\langle R_{4,0},\left(
\begin{array}{lll}
\mathcal{F}_{1}\\
\mathcal{F}_{2}
\end{array}
\right)\rangle+\langle R_{4,2},\left(
\begin{array}{lll}
\mathcal{F}_{5}\\
\mathcal{F}_{6}
\end{array}
\right)\rangle-4i\langle R_{4,1},\left(
\begin{array}{lll}
\Omega_{1}\\
\Omega_{2}
\end{array}
\right) \rangle-P_{4,1}^{(3)}G_{5}-P_{4,1}^{(4)}G_{6};$$
$$E_{2}:=P_{4,1}^{(3)}(G_{1}+G_{5})+P_{4,1}^{(4)}(G_{2}+G_{6})+P_{4,2}^{(1)}G_{3}+P_{4,2}^{(2)}G_{4};$$
$$E_{3}:=\langle R_{4,1},\left(
\begin{array}{lll}
\mathcal{F}_{3}\\
\mathcal{F}_{4}
\end{array}
\right)\rangle+4i\langle R_{4,1},\left(
\begin{array}{lll}
\Omega_{1}\\
\Omega_{2}
\end{array}
\right) \rangle; $$
$$
\begin{array}{lll}
\mathcal{F}_{1}&:=&-2\phi^{\lambda}\eta^{2}
P_{3,0}^{(4)}-i6\phi^{\lambda}\xi^{2}P_{3,0}^{(3)}-2\phi^{\lambda}\eta
R_{3,0}^{(2)}-i6\phi^{\lambda}\xi
R_{3,0}^{(1)}+i\phi^{\lambda}_{\lambda}\eta^{2}P_{2,0}^{(1)}\\
& &-2\phi^{\lambda}\xi\eta
P_{2,0}^{(2)}-2\xi\eta R_{2,0}^{(2)}+i\eta^{2}R_{2,0}^{(1)}
-3i\xi^{2}R_{2,0}^{(1)}-3i\phi^{\lambda}_{\lambda}\xi^{2}P_{2,0}^{(1)};
\end{array}
$$
$$
\begin{array}{lll}
\mathcal{F}_{2}&:=&-2\phi^{\lambda}\xi\eta
P_{3,0}^{(3)}-2i\phi^{\lambda}\xi\eta P_{3,0}^{(4)}
-2\phi^{\lambda}\eta R_{3,0}^{(1)}-2i\phi^{\lambda}\xi
R_{3,0}^{(2)}-2\phi^{\lambda}_{\lambda}\xi\eta P_{2,0}^{(1)}\\ &
&-2\xi\eta
R_{2,0}^{(1)}+3i\eta^{2}R_{2,0}^{(2)}-i\xi^{2}R_{2,0}^{(2)}
+3i\phi^{\lambda}\eta^{2}
P_{2,0}^{(2)}-i\phi^{\lambda}\xi^{2}P_{2,0}^{(2)};
\end{array}
$$
$$
\begin{array}{lll}
\mathcal{F}_{3}&:=&-2(\phi^{\lambda})^{2}\eta
P_{2,0}^{(2)}-2\phi^{\lambda}\eta R_{2,0}^{(2)}
+\frac{3i}{4}\xi\eta^{2}\\
& &-\frac{i3}{4}\xi^{3}-i6\phi^{\lambda}\xi
R_{2,0}^{(1)}-i6\phi^{\lambda}\phi^{\lambda}_{\lambda}\xi
P_{2,0}^{(1)};
\end{array}
$$
$$
\begin{array}{lll}
\mathcal{F}_{4}&:=&-2\phi^{\lambda}\phi^{\lambda}_{\lambda}\eta
P_{2,0}^{(1)}-2\phi^{\lambda}\eta
R_{2,0}^{(1)}-\frac{3}{4}\xi^{2}\eta\\
&
&-2iP_{2,0}^{(2)}(\phi^{\lambda})^{2}\xi+\frac{3}{4}\eta^{3}-2i\phi^{\lambda}\xi
R_{2,0}^{(2)};
\end{array}
$$
$$\mathcal{F}_{5}:=i\phi^{\lambda}\eta^{2}-i3\phi^{\lambda}\xi^{2};$$
$$\mathcal{F}_{6}:=-2\phi^{\lambda}\xi\eta;$$
$$\Omega_{1}:=3\phi^{\lambda}\xi R_{2,0}^{(1)}-i\phi^{\lambda}\eta R_{2,0}^{(2)},\ \Omega_{2}:=\phi^{\lambda}\xi R_{2,0}^{(2)}-i\phi^{\lambda}\eta R_{2,0}^{(1)};$$
$$
\begin{array}{lll}
G_{1}&:=&-\langle 2(\phi^{\lambda})^{2}\xi\eta
P_{0,2}^{(2)}\rangle-2\langle \phi^{\lambda}\xi\eta
R_{0,2}^{(2)}\rangle-\frac{3i}{4}\langle \xi^{2}\eta^{2}\rangle\\&
&+\langle 6i\phi^{\lambda}\phi^{\lambda}_{\lambda}\xi^{2}\rangle
P_{0,2}^{(1)} +\frac{3i}{4}\langle \xi^{4}\rangle+6i\langle
\phi^{\lambda}\xi^{2}R_{0,2}^{(1)}\rangle;
\end{array}
$$
$$
\begin{array}{lll}
G_{2}&:=&2iP_{0,2}^{(2)}\langle
(\phi^{\lambda})^{2}\xi\eta\rangle-\langle
2\phi^{\lambda}\phi^{\lambda}_{\lambda}\eta^{2}\rangle
P_{0,2}^{(1)}-2\langle
\phi^{\lambda}\eta^{2}R_{0,2}^{(1)}\rangle\\
& &+\frac{3}{4}\langle \eta^{4}\rangle-\frac{3}{4}\langle
\xi^{2}\eta^{2}\rangle+i\langle 2 \phi^{\lambda}\xi\eta
R^{(2)}_{0,2}\rangle;
\end{array}
$$
$$
G_{3} :=-i\langle
\phi^{\lambda}\phi^{\lambda}_{\lambda}\eta^{2}\rangle+3i\langle
\phi^{\lambda}\phi^{\lambda}_{\lambda}\xi^{2}\rangle;
$$
$$G_{4}:=-\langle (\phi^{\lambda})^{2}\xi\eta\rangle;$$
$$G_{5}:=12i\langle \phi^{\lambda}\xi^{2},R_{2,0}^{(1)}\rangle-4\langle \phi^{\lambda}\xi\eta, R_{2,0}^{(2)}\rangle;$$
$$G_{6}:=-4\langle \phi^{\lambda}\eta^{2},R_{2,0}^{(1)}\rangle+4i\langle \phi^{\lambda}\xi\eta,
R_{2,0}^{(2)}\rangle$$ where we use $\langle g_{1}\rangle$, $\langle
g_{1},g_{2}\rangle$ to denote $\int g_{1}(x)dx$, $\int
g_{1}\bar{g}_{2}dx$ respectively for any function $g_{1}$ and
$g_{2}$.
\subsection{Expansion of $E_{n},\ n=1,2,3$}
In this section we first relate the terms $R_{4,1}$, $R_{4,2}$, $P_{4,2}^{(k)},\ k=1,2,3,4,$ to the function $R_{4,0}$ so that $ReZ_{4,3}=Re\langle R_{4,0},\Psi\rangle$ for some function $\Psi.$ Then after some manipulation we relate the function $\Psi$ to $N_{4,0},$ hence to $R_{4,0}.$ Recall that
the relation between $R_{4,0}$ and $N_{4,0}$ in Lemma ~\ref{LM:R40}.

The main result is the following proposition whose proof is in the
later part of this section.
\begin{proposition}\label{LM:e1}
\begin{equation}\label{eq:e2}
Re[E_{2}]=Re[E_{4,0}+E_{4,1}+P_{4,1}^{(3)}Y_{1}+P_{4,1}^{(4)}Y_{2}],
\end{equation} and
\begin{equation}\label{eq:y1y2}
Re[P_{4,1}^{(3)}Y_{1}+P_{4,1}^{(4)}Y_{2}]=Re[\langle
R_{4,0}^{(1)},-6i \mathcal{H}^{(2)}_{2}\rangle+\langle
R_{4,0}^{(2)},-6i \mathcal{H}_{2}^{(1)}\rangle],
\end{equation}
\begin{equation}\label{eq:e1}
\begin{array}{lll}
Re[E_{1}+E_{4,0}]&=&Re\langle R_{4,0}^{(1)},
-8i \mathcal{H}_{1}^{(2)}-2i \mathcal{H}_{2}^{(2)}-2i \mathcal{H}_{3}^{(2)}\rangle\\
& &+Re\langle R_{4,0}^{(2)}, -8i \mathcal{H}^{(1)}_{1}-2i
\mathcal{H}_{2}^{(1)}-2i \mathcal{H}_{3}^{(1)}\rangle,
\end{array}
\end{equation}
\begin{equation}\label{eq:e3}
Re[E_{3}+E_{4,1}]=Re[\langle R_{4,0}^{(1)},-6i
\mathcal{H}_{3}^{(2)}\rangle+\langle
R_{4,0}^{(2)},-6i\mathcal{H}_{3}^{(1)}\rangle]
\end{equation} where the functions $Y_{1}$, $Y_{2}$, $E_{4,0}$ and
$E_{4,1}$ are defined below.
\end{proposition}

In the decomposition of $E_{2}$ above we used the following
constants
$$
E_{4,0}:=-4i\langle R_{4,0}^{(1)},
M_{1}^{(1)}P_{0,2}^{(2)}+M_{1}^{(2)}P_{0,2}^{(1)}\rangle-4i\langle
R_{4,0}^{(2)},
M_{2}^{(1)}P_{0,2}^{(2)}+M_{2}^{(2)}P_{0,2}^{(1)}\rangle$$
$$E_{4,1}:= -4i\langle R_{4,1}^{(1)},
M_{3}^{(1)}P_{0,2}^{(2)}+M_{3}^{(2)}P_{0,2}^{(1)}\rangle-4i\langle
R_{4,1}^{(2)},
M_{4}^{(1)}P_{0,2}^{(2)}+M_{4}^{(2)}P_{0,2}^{(1)}\rangle,
$$
$$
Y_{1}:=-4iW_{1}^{(1)}P_{2,0}^{(2)}-4iW_{1}^{(2)}P_{2,0}^{(1)}+G_{1}+G_{5},
$$ and $$
Y_{2}:=-4iW_{2}^{(1)}P_{2,0}^{(2)}-4iW_{2}^{(2)}P_{2,0}^{(1)}+G_{2}+G_{6}.
$$

Recall the notation
$\delta^{'}(\lambda):=\frac{1}{2}\frac{d}{d\lambda}\langle
\phi^{\lambda},\phi^{\lambda}\rangle.$ Before starting proving the
proposition, we compute the explicit form of $P_{4,2}^{(1)},
P_{4,2}^{(2)}:$ \begin{equation}\label{eq:p142}
\begin{array}{lll}
& &-2i\epsilon(\lambda)\delta^{'}(\lambda)P_{4,2}^{(1)}\\&=&\langle
R_{4,0}^{(1)}, \mathcal{M}_{1}^{(1)}\rangle + \langle
R_{4,0}^{(2)},\mathcal{M}_{2}^{(1)}\rangle+\langle
R_{4,1}^{(1)},\mathcal{M}_{3}^{(1)}\rangle +\langle
R_{4,1}^{(2)},\mathcal{M}_{4}^{(1)}\rangle\\
& &+P_{4,1}^{(3)}W_{1}^{(1)}+P_{4,1}^{(4)}W_{2}^{(1)}+iU.
\end{array}
\end{equation} where, $U$ is a real constant and $iU$ includes all the contributions from $F_{4}$ (see Lemma ~\ref{LM:junk}),
$$\mathcal{M}_{1}^{(1)}:=2(\phi^{\lambda})^{3}P_{2,0}^{(2)}+2(\phi^{\lambda})^{2}R_{2,0}^{(2)}-\frac{i}{2}\phi^{\lambda}\xi\eta$$
$$\mathcal{M}_{2}^{(1)}:=2(\phi^{\lambda})^{2}\phi^{\lambda}_{\lambda}P_{2,0}^{(1)}+2(\phi^{\lambda})^{2}R_{2,0}^{(1)}+\frac{1}{4}\phi^{\lambda}\xi^{2}-\frac{3}{4}\phi^{\lambda}\eta^{2}$$
$$\mathcal{M}_{3}^{(1)}:=-i(\phi^{\lambda})^{2}\eta\ \ \mathcal{M}_{4}^{(1)}:=(\phi^{\lambda})^{2}\xi,$$
$$W_{1}^{(1)}:=i\langle
(\phi^{\lambda})^{2}\xi\eta\rangle,\ \ W_{2}^{(1)}:=\langle
(\phi^{\lambda})^{2}\xi\eta\rangle,$$

\begin{equation}\label{eq:p242}
\begin{array}{lll}
-2i\epsilon(\lambda)\delta^{'}(\lambda)P_{4,2}^{(2)}&=&\langle
R^{(1)}_{4,0},\mathcal{M}^{(2)}_{1}+\frac{\mathcal{M}_{1}^{(1)}}{2i\epsilon(\lambda)}\rangle+\langle
R^{(2)}_{4,0},\mathcal{M}^{(2)}_{2}+\frac{{M}_{2}^{(1)}}{2i\epsilon(\lambda)}\rangle\\
& &+\langle
R^{(1)}_{4,1},\mathcal{M}^{(2)}_{3}+\frac{\mathcal{M}_{3}^{(1)}}{2i\epsilon(\lambda)}\rangle+\langle
R^{(2)}_{4,1},\mathcal{M}^{(2)}_{4}+\frac{\mathcal{M}_{4}^{(1)}}{2i\epsilon(\lambda)}\rangle\\
&
&+P_{4,1}^{(3)}[W^{(2)}_{1}-\frac{W^{(1)}_{1}}{2i\epsilon(\lambda)}]+P_{4,1}^{(4)}[W^{(2)}_{2}-\frac{W_{2}^{(1)}}{2i\epsilon(\lambda)}]+U,
\end{array}\end{equation}
where $U$ is a real constant including all the contributions from
$F_{3}$ (see Lemma ~\ref{LM:junk}), $$W^{(2)}_{1}:=-\langle
3\phi^{\lambda}\phi^{\lambda}_{\lambda}\xi^{2}\rangle,\
W^{(2)}_{2}:=-i\langle
\phi^{\lambda}\phi^{\lambda}_{\lambda}\eta^{2}\rangle,$$
$$
\mathcal{M}^{(2)}_{1}:=-6\phi^{\lambda}(\phi^{\lambda}_{\lambda})^{2}P_{2,0}^{(1)}-6\phi^{\lambda}\phi^{\lambda}_{\lambda}R_{2,0}^{(1)}-\frac{3}{4}\phi^{\lambda}_{\lambda}\xi^{2}+\frac{1}{4}\phi^{\lambda}_{\lambda}\eta^{2}.
$$
$$\mathcal{M}^{(2)}_{2}:=-2(\phi^{\lambda})^{2}\phi^{\lambda}_{\lambda}P_{2,0}^{(2)}-2\phi^{\lambda}\phi^{\lambda}_{\lambda}R_{2,0}^{(2)}+\frac{i}{2}\phi^{\lambda}_{\lambda}\xi\eta ,$$
$$\mathcal{M}^{(2)}_{3}:=-3\phi^{\lambda}\phi^{\lambda}_{\lambda}\xi,\ \mathcal{M}^{(2)}_{4}:=i\phi^{\lambda}\phi^{\lambda}_{\lambda}\eta,$$
In the following lemma we relate $Y_{1}$ and $Y_{2}$ to
$P_{3,0}^{(3)}$ and $P_{3,0}^{(4)}.$
\begin{lemma}
\begin{equation}\label{eq:y1y2form}
Y_{1}=6i\epsilon(\lambda)\langle
\xi,\eta\rangle[3iP_{3,0}^{(4)}-P_{3,0}^{(3)}],\
Y_{2}=6i\epsilon(\lambda)\langle \xi,\eta\rangle
[3iP_{3,0}^{(3)}+P_{3,0}^{(4)}].
\end{equation}
\end{lemma}
\begin{proof} Using Equations (~\ref{eq:z1}) (~\ref{eq:z2}) to
obtain
\begin{equation}\label{eq:p303}
-8\epsilon(\lambda)\langle \xi,\eta\rangle P_{3,0}^{(3)}=-3i\langle
ImN_{3,0},\eta\rangle-\langle ReN_{3,0},\xi\rangle
\end{equation} and
\begin{equation}\label{eq:p304}
8\epsilon(\lambda)\langle \xi,\eta\rangle P_{3,0}^{(4)}=-3i\langle
ReN_{3,0},\xi\rangle+\langle ImN_{3,0},\eta\rangle.
\end{equation}
On the other hand we observe $$
\begin{array}{lll}
Y_{1}&=&\frac{3}{4}i\langle \xi^{4}\rangle-\frac{3i}{4}\langle
\xi^{2}\eta^{2}\rangle+6\langle (\phi^{\lambda})^{2}\xi\eta\rangle
P_{2,0}^{(2)}+6\langle \phi^{\lambda}\xi\eta
R_{2,0}^{(2)}\rangle+18i\langle
\phi^{\lambda}\xi^{2}R_{2,0}^{(1)}\rangle\\
& &+18i\langle \phi^{\lambda}\phi^{\lambda}_{\lambda}\xi^{2}\rangle
P_{2,0}^{(1)}\\
&=&6i\langle ReN_{3,0},\xi\rangle
\end{array}
$$ and $$
\begin{array}{lll}
Y_{2}&=&\frac{3}{4}\langle \eta^{4}\rangle-\frac{3}{4}\langle
\xi^{2}\eta^{2}\rangle-6iP_{2,0}^{(2)}\langle
(\phi^{\lambda})^{2}\xi\eta\rangle-6\langle
\phi^{\lambda}\eta^{2}R_{2,0}^{(1)}\rangle-6i\langle
\phi^{\lambda}\xi\eta R_{2,0}^{(2)}\rangle\\
& &-6\langle
\phi^{\lambda}\phi^{\lambda}\eta^{2}\rangle P_{2,0}^{(1)}\\
&=&-6i\langle ImN_{3,0},\eta\rangle.
\end{array}
$$ These observations together with (~\ref{eq:p303}) and (~\ref{eq:p304}) imply (~\ref{eq:y1y2form}).
\end{proof}
\textbf{Proof of Equations (~\ref{eq:e2}) and (~\ref{eq:y1y2})} By
the definition of $E_{2}$ it is sufficient to prove (~\ref{eq:e2})
by proving
\begin{equation}\label{eq:e2part}
Re[P_{4,2}^{(1)}G_{3}+P_{4,2}^{(2)}G_{4}]=Re[E_{4,0}+E_{4,1}+P_{4,1}^{(3)}(Y_{1}-G_{1}-G_{5})+P_{4,1}^{(4)}(Y_{2}-G_{2}-G_{6})].
\end{equation}
First we observe that
$$G_{3}-\frac{1}{2i\epsilon(\lambda)}G_{4}=-8\epsilon(\lambda)\delta^{'}(\lambda) P_{2,0}^{(2)},\ \ G_{4}=-8\epsilon(\lambda)\delta^{'}(\lambda)P_{2,0}^{(1)}.$$
Lemma ~\ref{LM:keyfact} implies that $P^{(1)}_{2,0}$ and
$P^{(2)}_{2,0}$ are real and purely imaginary respectively, this
fact together with Equations (~\ref{eq:p142}) (~\ref{eq:p242})
implies (~\ref{eq:e2part}). The computation is straightforward, but
tedious, hence we omit the detail.

Now we turn to (~\ref{eq:y1y2}). We start with analyzing
$P_{4,1}^{(3)}$ and $P_{4,1}^{(4)}$. Equations (~\ref{eq:z1})
(~\ref{eq:z2}) and the properties of the equation for $\dot{z}$ in
(~\ref{FGRZ}) imply
\begin{equation}\label{eq:p413}
-8\epsilon(\lambda)\langle \xi,\eta\rangle P_{4,1}^{(3)}=-3i\langle
ImN_{4,1},\eta\rangle-\langle ReN_{4,1},\xi\rangle+U
\end{equation} and
\begin{equation}
8\epsilon(\lambda)\langle \xi,\eta\rangle P_{4,1}^{(4)}=-3i\langle
ReN_{4,1},\xi\rangle+\langle ImN_{4,1},\eta\rangle+iU
\end{equation} where $U$ stands for real (different) constants including all the contribution from the terms $F_{1}$ and $F_{2}$ (see Lemma ~\ref{LM:junk}).
Moreover we compute to get
\begin{equation}\label{eq:n41}
N_{4,1}=\left(
\begin{array}{lll}
ImN_{4,1}\\
-ReN_{4,1}
\end{array}
\right)=\left(
\begin{array}{lll}
i\phi^{\lambda}\eta R_{4,0}^{(1)}+\phi^{\lambda}\xi R_{4,0}^{(2)}\\
-3\phi^{\lambda}\xi R_{4,0}^{(1)}-i\phi^{\lambda}\eta R_{4,0}^{(2)}
\end{array}
\right)+U_{1}
\end{equation} where $U_{1}$ is a vector function such that $iU_{1}$
is admissible. Thus we have $$\langle ReN_{4,1},\xi\rangle=\langle
R_{4,0},\left(
\begin{array}{lll}
3\phi^{\lambda}\xi^{2}\\
-i\phi^{\lambda}\xi\eta
\end{array}
\right)\rangle+U $$ and $$\langle ImN_{4,1},\eta\rangle=\langle
R_{4,0},\left(
\begin{array}{lll}
-i\phi^{\lambda}\eta^{2}\\
\phi^{\lambda}\xi\eta
\end{array}
\right)\rangle+iU$$

Collecting the computation above and using the properties of $Y_{1}$
and $Y_{2}$ in (~\ref{eq:y1y2form}) and the facts that
$\overline{P_{3,0}^{(3)}}=P_{3,0}^{(3)},$
$\overline{P_{3,0}^{(4)}}=-P_{3,0}^{(4)}$ in Lemma ~\ref{LM:keyfact}
we have (~\ref{eq:y1y2}).
\begin{flushright}
$\square$
\end{flushright}
\textbf{Proof of Equation (~\ref{eq:e1})} We start with rewriting
$\langle R_{4,2},\left(
\begin{array}{lll}
\mathcal{F}_{5}\\
\mathcal{F}_{6}
\end{array}
\right)\rangle$ to relate the function $\left(
\begin{array}{lll}
\mathcal{F}_{5}\\
\mathcal{F}_{6}
\end{array}
\right)$ to $R_{2,0}.$ By the definition of $R_{4,2}$ in Lemma
~\ref{LM:keyfact} we have
$$R_{4,2}=-(L(\lambda)+2i\epsilon(\lambda))^{-1}P_{c}[N_{4,2}+U_{1}]$$ where $U_{1}$ is a vector function such that $iU_{1}$ is
admissible and includes all the contributions of $\mathcal{G}$ in
(~\ref{RAfProj}) (see Lemma ~\ref{LM:junk} for $N=3$).

Recall that $\sigma_{1}=\left(
\begin{array}{lll}
0&-1\\
1&0
\end{array}
\right)$. The facts
$\sigma_{1}(L(\lambda))^{*}\sigma_{1}=L(\lambda)$ and
$\sigma^{2}_{1}=-1$ imply $$ \langle R_{4,2},\left(
\begin{array}{lll}
\mathcal{F}_{5}\\
\mathcal{F}_{6}
\end{array}\right)\rangle
=\langle N_{4,2}
,-\sigma_{1}(L(\lambda)+2i\epsilon(\lambda))^{-1}P_{c}\sigma_{1}\left(
\begin{array}{lll}
\mathcal{F}_{5}\\
\mathcal{F}_{6}
\end{array}
\right)\rangle+iU
$$ where $U$ is a real constant generated by $U_{1}$ above (see Lemma ~\ref{LM:admissibility1}). Moreover we observe that
$$
\begin{array}{lll}
R_{2,0}&=&-(L(\lambda)+2i\epsilon(\lambda))^{-1}P_{c}N_{2,0}\\
&=&\frac{i}{4}(L(\lambda)+2i\epsilon(\lambda))^{-1}P_{c}\sigma_{1}\left(
\begin{array}{lll}
\mathcal{F}_{5}\\
\mathcal{F}_{6}
\end{array}
\right)
\end{array}
$$ consequently $$\langle R_{4,2},\left(
\begin{array}{lll}
\mathcal{F}_{5}\\
\mathcal{F}_{6}
\end{array}
\right)\rangle=-4i\langle N_{4,2}, \sigma_{1}R_{2,0}\rangle+iU.$$

In the next we decompose $N_{4,2}$ into various terms. By direct
computation we have
$$
\begin{array}{lll}
N_{4,2}&=&\mathcal{M} R_{4,0} +\left(
\begin{array}{lll}
i\phi^{\lambda}\eta&\phi^{\lambda}\xi\\
-3\phi^{\lambda}\xi&-i\phi^{\lambda}\eta
\end{array}
\right)R_{4,1}\\
& &+\left(
\begin{array}{lll}
i\phi^{\lambda}\xi\eta&\phi^{\lambda}\xi\eta \\
-3\phi^{\lambda}\xi^{2}&-i\phi^{\lambda}\eta^{2}
\end{array}
\right)\left(
\begin{array}{lll}
P_{4,1}^{(3)}\\
P_{4,1}^{(4)}
\end{array}
\right)+U_{1}
\end{array}
$$ where $U_{1}$ is a vector function with $iU_{1}$ being admissible,
$\mathcal{M}$ is an $2\times 2$ matrix defined as
$$\left(
\begin{array}{lll}
2\phi^{\lambda}R_{0,2}^{(2)}+\frac{i}{2}\xi\eta+2(\phi^{\lambda})^{2}P_{0,2}^{(2)}, &  \frac{1}{4}\xi^{2}-\frac{3}{4}\eta^{2}+2\phi^{\lambda}\phi^{\lambda}_{\lambda}P_{0,2}^{(1)}+2\phi^{\lambda}R_{0,2}^{(1)}\\
\frac{3}{4}\xi^{2}+6\phi^{\lambda}\phi^{\lambda}_{\lambda}P_{0,2}^{(1)}-\frac{1}{4}\eta^{2}+6\phi^{\lambda}R_{0,2}^{(1)},&
2\phi^{\lambda}R_{0,2}^{(2)}+2(\phi^{\lambda})^{2}P_{0,2}^{(2)}+\frac{i}{2}\xi\eta
\end{array}
\right).$$

The definitions of $E_{1}$ after (~\ref{eq:expZ43}) yields
$$E_{1}=\langle R_{4,0},\left(
\begin{array}{lll}
\mathcal{F}_{1}\\
\mathcal{F}_{2}
\end{array}
\right)\rangle-4i\langle \mathcal{M}R_{4,0},
\sigma_{1}R_{2,0}\rangle+iU$$ which together with the definition of
$E_{4,0}$, Lemma ~\ref{LM:keyfact} and the fact $R_{m,n},\ m,n\leq
3$ is admissible implies Equation (~\ref{eq:e1}). The computation is
straightforward, hence we omit the detail.
\begin{flushright}
$\square$
\end{flushright}
\textbf{Proof of Equation (~\ref{eq:e3})} The definitions of $E_{3}$
and $E_{4,1}$ imply
$$E_{3}+E_{4,1}=\langle R_{4,1}, K
 \rangle+iU$$ with $U$ being a real constant and $$K:=\left(
\begin{array}{lll}
\mathcal{F}_{3}-4i\Omega_{1}+4i[M^{(1)}_{3}P_{0,2}^{(2)}+M_{3}^{(2)}P_{0,2}^{(1)}]\\
\mathcal{F}_{4}-4i\Omega_{2}+4i[M_{4}^{(1)}P_{0,2}^{(2)}+M_{4}^{(2)}P_{0,2}^{(1)}]
\end{array}
\right).$$ In the next we relate the function $K$ to $R_{3,0}.$

First we have
$$R_{4,1}=-(L(\lambda)+3i\epsilon(\lambda))^{-1}P_{c}[N_{4,1}+U_{1}]$$ where as usual $U_{1}$ is a vector function such
that $iU_{1}$ is admissible and includes all the contributions from
$\mathcal{G}_{1}$ from (~\ref{RAfProj}). On the other hand by Lemma
~\ref{LM:keyfact} and direct computation we have
$-6i\sigma_{1}N_{3,0}=K $, thus $$
\begin{array}{lll}
E_{3}+E_{4,1}&=&6i\langle R_{4,1}, \sigma_{1}N_{3,0}\rangle\\
&=&6i\langle
N_{4,1},-(L(\lambda)^{*}+3i\epsilon(\lambda))^{-1}P_{c}^{*}\sigma_{1}N_{3,0}\rangle+iU
\end{array}
$$ with $U$ being a real constant generated by $U_{1}$ above (see Lemma ~\ref{LM:admissibility1}). By the facts that
$\sigma_{1}(L(\lambda))^{*}\sigma_{1}=L(\lambda)$,
$\sigma^{2}_{1}=-1$, $-\sigma_{1}P_{c}^{*}\sigma_{1}=P_{c} $ and the
definition
$R_{3,0}=-(L(\lambda)+3i\epsilon(\lambda))^{-1}P_{c}N_{3,0}$ we have
$$E_{3}+E_{41}=-6i\langle N_{4,1},\sigma_{1}R_{3,0}\rangle.$$ By the
explicit formula for $N_{4,1}$ in (~\ref{eq:n41}), we obtain
$$
\begin{array}{lll}
E_{3}+E_{41}=\langle R_{4,0},-6i\left(
\begin{array}{lll}
\mathcal{H}_{3}^{(2)}\\
\mathcal{H}_{3}^{(1)}
\end{array}
\right)\rangle+iU
\end{array}
$$ for some real constant $U.$ The proof is complete.
\begin{flushright}
$\square$
\end{flushright}
\subsection{Proof of Theorem ~\ref{THM:mainresult}}
\begin{proof}
Recall the explicit form of $Z_{4,3}$ and $R_{4,0}$ in
(~\ref{eq:expZ43}). By Equations (~\ref{eq:e2})-(~\ref{eq:e3}) and
the observation that $N_{4,0}=\left(
\begin{array}{lll}
\displaystyle\sum_{n=1}^{3}\mathcal{H}_{n}^{(1)},\
-\displaystyle\sum_{n=1}^{3}\mathcal{H}_{n}^{(2)}
\end{array}
\right)^{T}$ in (~\ref{eq:ImN40}) we have
$$
\begin{array}{lll}
ReZ_{4,3}&=&\frac{8}{\langle \xi,\eta\rangle}Im\langle
\sigma_{1}(L(\lambda)+4i\epsilon(\lambda)+0)^{-1}N_{4,0},N_{4,0}\rangle\\
&=&\frac{8}{\langle \xi,\eta\rangle}Im\langle
\sigma_{1}(L(\lambda)-4i\epsilon(\lambda)-0)^{-1}N_{4,0},N_{4,0}\rangle
\end{array}
$$ to complete the proof.
\end{proof}

\end{document}